\documentclass[onefignum,onetabnum]{siamart190516}

\usepackage{geometry}
\usepackage{amssymb}
\usepackage[utf8]{inputenc}
\usepackage{graphicx} 
\usepackage{subcaption}
\usepackage{tikz,pgfplots}
\usetikzlibrary{patterns}
\usetikzlibrary{arrows, arrows.meta}

\usepackage[colorinlistoftodos]{todonotes}
\usepackage[capitalize,nameinlink]{cleveref}

\graphicspath{{figures/}}
\DeclareGraphicsExtensions{.pdf,.eps,.png,.jpg,.jpeg}

\title{Barrier Simulations and Experimental Calculations Using Cell Merging Method}
\author{Chanyang Ryoo \and Kyle T. Mandli}
\usepackage[colorinlistoftodos]{todonotes}
\usepackage[capitalize,nameinlink]{cleveref}

\ifpdf
\hypersetup{
  pdftitle={Application of },
  pdfauthor={C. Ryoo, and K. T. Mandli}
}
\fi

\begin{document}

\maketitle

\begin{abstract}
   We apply the cell merging method to a model shallow water problem with a permeable boundary. We use a cut cell approach which is more easily and systematically scalable with different shapes of boundaries. The novel cell merging method presented in this paper uses both wave propagation algorithm and gradient reconstruction for second order corrections, along with minmod and Barth-Jespersen limiters. We observe second order convergence in two model problems, one with a linear boundary and other with a composite, V-shaped boundary. We assess the effectiveness of these two boundaries by doing a realistic scenario including an island and observing inundation at the peak of the island.
\end{abstract}





\section{Introduction}
Barrier simulations are difficult due to the barriers' inherently thin structure \cite{o2020hydrodynamic,lindemer2010numerical}. This means that more resolution is required at the barriers. We propose a model problem that approximates a barrier using a line interface embedded on a Cartesian grid. Instead of using unstructured meshing, we use this embedded grid framework as the latter removes the need of generating new meshes for new barriers and allows for straightforward scaling and generalization.

This framework requires the use of cut cell methods \cite{BERGER20171} because of the embedded boundary cutting through some of the grid cells. There are many cut cell methods available. All of them handle the CFL restriction at the cut cells. The benefits of using cut cell methods over resolving the barrier are increased time step sizes and reduced number of steps to a desired final $t$. Here we employ the \emph{cell merging} method as it is relatively easy to implement and conceptually accessible, unlike other methods which are not so theoretically straightforward such as the state redistribution \cite{BERGER2020109820}. Also, as will be shown, applying a second order method will be relatively simple to implement.

The equations we solve are the shallow water equations, used in tsunami, storm, flooding simulations \cite{BERGER20111195,jmse8090725}. We test our method on two model problems: one with a linear barrier and other with a piecewise linear barrier that forms the shape $V$ at an obtuse angle. As a side note, these shapes are the proposed shapes of the NYC surge barrier that will connect Rockaways and Sandy Hook. The barriers have height parameter $\beta$ which will determine whether an incoming wave can overtop them. We show both complete blockage and overtopping examples.

\section{Model equation and problem} \label{sec1}
We present the underlying partial differential equations that we solve numerically and the grid setup of the barrier problems with some of their specifications. We will show the numerical results of these problems in \cref{sec4}.
\subsection{2D shallow water equations}
In many cases shallow water equations (SWE) have been used to simulate flooding, storms, and tsunamis \cite{martin2010lake,BERGER20111195,zhang2013transition}. They are a system of hyperbolic partial differential equations that is derived by depth-averaging the Navier-Stokes equation. There are ways to add source terms to construct a geophysical version of SWE to simulate wave motion on a rotating sphere like the earth \cite{calhoun2008logically}, but we focus on the simple SWE with varying bathymetry being our only source terms.

The shallow water equations are a simplification of Navier-Stokes but are more readily understood as a set of conservation equations \cite{GDavid}:
\begin{align}
& h_t + (hu)_x + (hv)_y = 0 \\
& (hu)_t + (\frac{1}{2}gh^2 + hu^2)_x + (huv)_y = -ghb_x\\
& (hv)_t + (huv)_x + (\frac{1}{2}gh^2+hv^2)_y = -ghb_y,
\end{align}
for $(x,y) \in \Omega \subset \mathbb{R}^2$, where $h$ represents the height of water, $u$ and $v$ the $x$-velocity and $y$-velocity, $g$ the gravitational constant, and $b$ the bathymetry. In the first equation, we have conservation of mass, where the time derivative of the height (or mass), $h_t$, over a certain region $\Omega$ only depends on the mass flux $(hu)_x$ and $(hu)_y$ through $\Omega$. In the second and third equations, if we have no bathymetric variation ($\nabla b \equiv 0$), we have conservation of momentum, where the time derivative of momentum $(hu)_t$, $(hv)_t$ over $\Omega$ only depends on the momentum flux $(1/2 gh^2 + hu^2)_x + (huv)_y$, $(huv)_x + (1/2gh^2 + hv^2)_y$ through $\Omega$, respectively. In the presence of bathymetric variation (-$ghb_x$ or -$ghb_y$) there is loss of conservation of momentum. We also note that the bathymetry $b$ is a time independent variable.

We can thus represent the SWE in the form of a conservative equation \cref{conseqn} as derived in \cite{leveque2002finite}. We use this form for our numerical method later in \cref{sec2}. Let $q$ represent the state vector, $[h,hu,hv]$, and $f(q)$ the (mass and momentum) flux vector in $x$-direction, $[q_2, \frac{1}{2}gq_1^2 + \frac{q_2^2}{q_1},\frac{q_2q_3}{q_1}]$, and $g(q)$ the flux vector in $y$-direction, $[q_3,\frac{q_2q_3}{q_1},\frac{1}{2}gq_1^2 + \frac{q_3^2}{q_1}]$, and $\Psi(q,b)$ the bathymetric source term $ [0,-gq_1b_x, -gq_1b_y]$. Then we have for our SWE,
\begin{align}
    q_t + f(q)_x + g(q)_y = \Psi(q,b).
    \label{conseqn}
\end{align}

\subsection{Barrier representation on Cartesian grid}
\begin{figure}[h!]
\begin{subfigure}[b]{0.5\textwidth}
\centering
\begin{tikzpicture}
\draw[step=0.4cm] (0,0) grid (4,4);
\node at (-0.1, 1.2) {$1$};
\draw[red,thick] (0,1.2)--(4,2.5);
\node at (4.1, 2.5) {$2$};
\end{tikzpicture}
\caption{Model problem 1: slanted barrier.}
\label{fig:2d_setup1}
\end{subfigure}
 \hfill
\begin{subfigure}[b]{0.5\textwidth}
\centering
\begin{tikzpicture}
\draw[step=0.4cm] (0,0) grid (4,4);
\node at (-0.1, 2.7) {$1$};
\draw[red,thick] (0,2.7)--(2,1.27)--(4,2.7);
\node at (2, 1.2) {$2$};
\node at (4.1, 2.7) {$3$};

\end{tikzpicture}
\caption{Model problem 2: $V$-barrier}
\label{fig:2d_setup2}
\end{subfigure}
\caption{Grid setup of two model barrier problems. $1$, $2$, and $3$ denote the vertices of the barriers.}
\label{fig:2d_setup}
\end{figure}
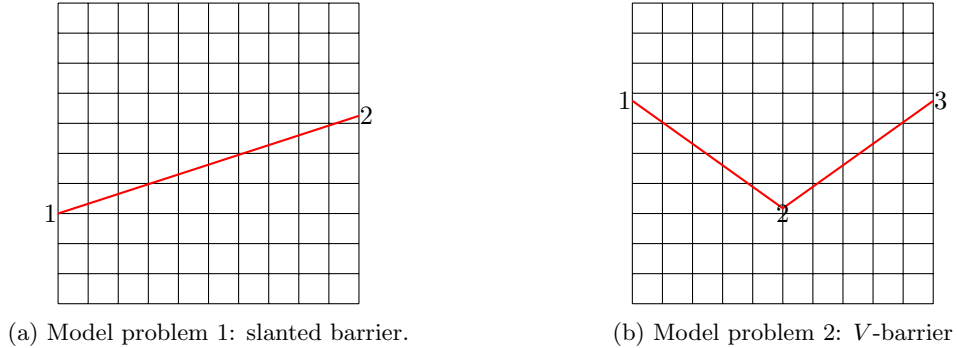

We now present the two model problems by showing their grid setup in \cref{fig:2d_setup}. The grid is intentionally coarse to show the nature of the cut cells and the small cell problem they pose.

In red we have our representation of the barrier as line interfaces. Note that the cut produces two states on either side of the line in the same cell position. This has both research and coding implications. It has research implications as this is the only instance known by the authors where there are state cells on both sides of a cut cell. In all cases that deal with cut cells \cite{CAUSON2000545,TUCKER2000591,may2017explicit}, one side of the cut is an impermeable solid where there is no state variable and the other half is the state cell. This also has coding implications as we must have two arrays to keep track of the upper state values and the lower state values in the same grid cell position, whereas in the one-state only problems, there is only one grid with single state values everywhere (except the solid non-state cells).

Each grid cell has four values: $Q=[H, HU, HV]$ and $B$ (capitalized because they represent volume averages), for which we can set initial conditions. The only exception is that near the barrier, we restrict the bathymetry level $B$ to be uniform throughout the cut cells and their direct neighboring cells. This is to keep conservation while using the CM method, as the method is developed \cite{BERGER2020109820} for hyperbolic equations without any source term (e.g. bathymetric variation). However, this is a reasonable limitation as in reality, surge barriers are most likely built upon a level seabed. Also, since the height of the barrier above the water surface level should be same all along the barrier, we avoid the complexity of having to set different barrier heights $\beta$ on different cut cells.

\subsubsection{Computational aspects}
There are several input parameters for our model problem. First, we have the vertices of the barriers $1$, $2$, and $3$ (\cref{fig:2d_setup}) which determine the location and geometry of the linear and $V$-shaped barrier. Then we have the barrier height parameter $\beta$, which is the height measured from the bathymetry level of the cut cells to the top of the barrier. For our model problems, we set vertex $1 = (0,0.3)$ and vertex $2=(0,0.653)$ for the linear barrier case ($20^\circ$) and vertex $1=(0,0.72)$, vertex $2=(0.5,0.412)$ and vertex $3=(1,0.72)$ for the V barrier case ($117^\circ$). As we will see in the numerical examples, we alternate $\beta$ between $1.5$ and $5$ to observe both overtopping and reflection.

Note also both the arbitrary locations and the arbitrarily small size of the cut cells that are produced by the barriers (\cref{fig:2d_setup}). We compute the locations of cut cells by using the vertices $A$, $B$, and $C$ of the barrier segments and the grid mesh size $\Delta x=\Delta y$ to identify intersections between the barrier and the grid. We compute the area of a cut cell $\mathcal{A}$ by identifying its vertices $\{x_i,y_i\}_{i=1}^n$ and using the Shoelace formula, given by
\begin{align}
    \mathcal{A} = \frac{1}{2} \big | \big ( \sum_{i=1}^{n-1} x_iy_{i+1} \big ) + x_ny_1 - \big (\sum_{i=1}^{n-1} x_{i+1}y_i \big ) - x_1y_n \big |.
\end{align}
We compute the cut cell indices and areas in a preprocessing step given the inputs of vertices $A$, $B$, and $C$, and use the computed data during the numerical solving of equations to avoid having to repeatedly perform these calculations in solution runtime.

\subsubsection{Second Order Information}
To perform second order methods, we also need to compute gradients $\nabla Q_{i,j}$ to do linear approximations. We use the gradients to approximate the solution values at cell edge midpoints at the barrier. In order to compute $\nabla Q_{i,j}$ we will use the least square gradient reconstruction method (LSQ). To use LSQ, we compute (1) each midpoint of the cut cell grid edge $(x_e,y_e)$, (2) centroid of the cut cell $(x_i,y_i)$, and (3) distance vectors $\Delta \mathbf{r}_{\{N_{i,j}\}}$ to their neighboring cell centroids (five-point stencil in \cref{fig:stencil_grad}). We pre-compute these information for each grid, and use them to do least square fitting on the following equation:
\begin{align}
    Q_{\{N_{i,j}\}} - Q_{i,j} = \Delta \mathbf{r}_{\{N_{i,j}\}} \nabla Q_{i,j} ,
\end{align}
where $\{N_{i,j}\}$ makes up the stencil for $(i,j)$ with $|\{N_{i,j}\ : N_1, ..., N_n\}|=n$ ,
\begin{align*}
    Q_{\{N_{i,j}\}} - Q_{i,j} =\begin{bmatrix} Q_{N_i} - Q_{i,j} \\ \vdots \\ Q_{N_n}-Q_{i,j} \end{bmatrix}
\end{align*}
and dimension thus being $n \times 3$,
and
\begin{align*}
    \Delta \mathbf{r}_{\{N_{i,j}\}} = \begin{bmatrix} (x_{N_i}-x_i), (y_{N_i}-y_j) \\ \vdots \\ (x_{N_n}-x_i), (y_{N_n}-y_j) \end{bmatrix}
\end{align*}
with dimension $n\times 2$, the entries $(x_i,y_j)$ being the centroid of cell $(i,j)$ and $(x_{N_i},y_{N_i})$ being the centroid of cell neighboring $N_{i}$, and finally
\begin{align*}
    \nabla Q_{i,j} = \begin{bmatrix} h_x , (hu)_x, (hv)_x \\ h_y, (hu)_y, (hv)y \end{bmatrix}_{(i,j)},
\end{align*}
being the $2\times 3$ matrix form of the approximation of the gradient at $(i,j)$.

With the gradient $\nabla Q_{i,j}$ computed we can `walk' from the centroid $(x_i,y_j)$ to the cell grid edges $(x_e,y_e)$ to get $Q^e_{i,j} = Q_{i,j} + \nabla Q_{i,j}^{T} (x_e-x_i,y_e-y_i).$ We note that this gradient calculation is used only at the barrier edge and elsewhere we use second order wave propagation method, which we discuss later. Also in \cref{cut_2nd} we show in more detail the stencil used to approximate gradients and the limiters used.

\section{Base numerical method: wave propagation} \label{sec2}
For our base method, we employ the wave propagation algorithm developed in \cite{LEVEQUE1997327}. We discuss its first order formulation and then its second.
\subsection{First order} \label{wp_1st}
We employ a finite volume method to solve the shallow water equations called wave propagation, or flux-based wave decomposition \cite{bale2003wave}. Instead of using fluxes at the edges of a cell, the method linearly decomposes the difference of flux vectors from either side of each edge. The matrices used for the linear decomposition are given by
\begin{align}
    A = \begin{bmatrix}
    1 & 0 & 1 \\
    U^*-\sqrt{gH^*} & 0 & U^* + \sqrt{gH^*} \\
    V^* & 1 & V^*
    \end{bmatrix},
    B = \begin{bmatrix}
    1 & 0 & 1 \\
    U^* & 1 & U^*  \\
    V^*-\sqrt{gH^*} & 0 & V^*+ \sqrt{gH^*}
    \end{bmatrix},
    \label{matrices}
\end{align}
where matrix $A$ is used to decompose difference of the flux vectors $f(q)$, and $B$ is used for $g(q)$, and the asterisked values represent special averages between $Q_{i,j}$ and its neighboring cells, which can either be Roe average or Einfeldt average \cite{einfeldt1988godunov,roe1981approximate}. Then the decomposition equations for the right and top edge become
\begin{align}
    f(Q_{i+1,j})-f(Q_{i,j})=A\gamma \label{diffcoeff}\\
    g(Q_{i,j+1})-g(Q_{i,j})=B\delta,\label{diffcoeff2}
\end{align}
where $\gamma$ and $\delta$ are coefficient vectors in $\mathbb{R}^3.$ The reason for using the wave propagation method is that it is a versatile method that can also be applied to nonconservative hyperbolic equations \cite{leveque2002finite} and also an easily implemented method as long as the decomposing matrices are known and well-behaved.
The reason for using these matrices $A$ and $B$ \cref{matrices} is that their columns are the eigenvectors of the Jacobians $f'(q)$ and $g'(q)$, respectively. The eigenvalues corresponding to the $p^{\text{th}}$ column vector of $A$ are $\{(\sigma_A)_p\}_{p=1}^{3} =\{U^*-\sqrt{gH^*}, U^*, U^*+\sqrt{gH^*}\}$, and the eigenvalues  corresponding to the $p^{\text{th}}$ column vector of $B$ are $\{(\sigma_B)_p\}_{p=1}^3 =\{V^*-\sqrt{gH^*},V^*, V^*+\sqrt{gH^*}\}$.

\subsubsection{Update formula using wave propagation}
Let $\{(\sigma_A)_p, A_{(\sigma_A)_p}\}$ and $\{(\sigma_B)_p,B_{(\sigma_B)_p}\}$ represent the $p^{\text{th}}$ eigenvalue-eigenvector pair of the Jacobians $A$ and $B$. Then we have the \emph{$A$-left going waves} and \emph{$A$-right going waves} defined by
\begin{align}
    A^-\Delta Q_{i+1/2,j}=\sum_{p : (\sigma_A)_p <0} \gamma_p (\sigma_A)_p  A_{(\sigma_A)_p}\\
    A^+\Delta Q_{i-1/2,j} = \sum_{p : (\sigma_A)_p >0} \gamma_p (\sigma_A)_p A_{(\sigma_A)_p},
\end{align}
where $\gamma$ is the difference coefficient vector in \cref{diffcoeff}.

The \emph{$B$-up} and \emph{$B$-down going waves} are given by:
\begin{align}
B^-\Delta Q_{i,j-1/2} = \sum_{p : (\sigma_B)_p < 0} \delta_p (\sigma_B)_p B_{(\sigma_B)_p} \\
B^+\Delta Q_{i,j+1/2} = \sum_{p: (\sigma_B)_p >0} \delta_p (\sigma_B)_p B_{(\sigma_B)_p},
\end{align}
where $\delta$ is the difference coefficient vector in \cref{diffcoeff2}.
Then the state update formula via the wave propagation method becomes:
\begin{align}
\label{waveupd}
    Q^{n+1}_{i,j} = Q^n_{i,j} - \frac{\Delta t}{\Delta x}(A^- \Delta Q_{i+1/2,j} + A^+ \Delta Q_{i-1/2,j}) - \frac{\Delta t}{\Delta y}(B^- \Delta Q_{i,j+1/2} + B^+ \Delta Q_{i,j-1/2}).
\end{align}

\subsubsection{Relationship between flux differencing form and wave propagation}
Although this wave propagation algorithm is different from standard flux differencing method, there is an equivalence between the two for hyperbolic conservation equations such as the SWE. The relation \cite{leveque2002finite} between numerically approximated flux $F_{i-1/2}$ and the waves is given by: $F_{i-1/2} = f(q_{i-1})+\mathcal{A}^-\Delta Q_{i-1/2}, \, F_{i-1/2} = f(q_i) - \mathcal{A}^+\Delta Q_{i-1/2} $
which gives us:
\begin{align*}
    F_{i+1/2} - F_{i-1/2} = \mathcal{A}^-\Delta Q_{i+1/2} + \mathcal{A}^+ \Delta Q_{i-1/2}.
\end{align*}

With this relation we can go back and forth from the flux difference form of the numerical update to the wave propagation form akin to what is done in \cite{berger2003h}. We use the flux difference form to do conservation calculations later on in our $h$-box methods and use the wave propagation method for wave redistribution and implementation of examples. Also note that the terminology of ``flux" is used for flux difference form but ``fluctuation" is used for the wave propagation algorithm, but we will use them interchangeably.

\subsection{Second order}\label{wp_2nd}
Second order methods are made from using a linear approximation on each cell. From the first order method, however, we can easily move towards a second order method by careful use of the speeds (eigenvalues) and the eigenbasis vectors.  The wave propagation algorithm provides an equivalent way of getting the higher order correction terms necessary for a second order numerical method without approximating solution slopes \cite{LEVEQUE1997327,george2008augmented}.

This is achieved by adding to the first order fluctuations a correction term $\tilde{F}_{i\pm 1/2,j}, \tilde{G}_{i,j\pm 1/2}$:
\begin{align} \label{eq:sec_ord_upd}
    Q_{i,j}^{n+1} = Q_{i,j}^n - \frac{\Delta t}{\Delta x}(A^-\Delta Q_{i+1/2,j} + A^+ \Delta Q_{i-1/2,j}) - \frac{\Delta t}{\Delta x}(\tilde{F}_{i+1/2,j}-\tilde{F}_{i-1/2,j})  \nonumber\\
    - \frac{\Delta t}{\Delta y}(B^-\Delta Q_{i,j+1/2} + B^+\Delta Q_{i,j-1/2}) - \frac{\Delta t}{\Delta y}(\tilde{G}_{i,j+1/2} - \tilde{G}_{i,j-1/2}).
\end{align}

The correction terms are then given by the following:
\begin{align}
    \tilde{F}_{i-1/2,j} = \frac{1}{2} \sum_{p=1}^m |(\sigma_A)_p|\big (1 - \frac{\Delta t}{\Delta x} |(\sigma_A)_p| \big ) \tilde{\gamma}_p A_{(\sigma_A)_p}, \\
    \tilde{G}_{i,j-1/2} = \frac{1}{2} \sum_{p=1}^m |(\sigma_B)_p|\big ( 1- \frac{\Delta t}{\Delta x} |(\sigma_B)_p| \big ) \tilde{\delta}_p B_{(\sigma_B)_p},
\end{align}
where $m$ is the dimension of the eigenspectrum, and the tilde over $\gamma_p, \delta_p$ is a limited version of the decomposition coefficient which we discuss in the next subsection. In a linear hyperbolic system $q_t+Aq_x =0,$ this formulation \cref{eq:sec_ord_upd} (without the limiting over the coefficients) can be shown to be equal to the Lax-Wendroff method in each dimension \cite{leveque2002finite}.

\subsubsection{Wave limiting}\label{wave_lim}
The tilde over the decomposition coefficients $\tilde{\gamma},\tilde{\delta}$ denote limiting of those coefficients. Instead of limiting flux vectors, these work by limiting the eigenbasis vectors composing the flux vector differences from the Riemann problem. The limiting of the coefficients are performed by comparing the the magnitude of $p^{th}$ eigenvectors arising from neighboring Riemann problems and applying the typical limiting functions to the ratio of magnitudes.

In symbols, we have the following: for a given Riemann problem between $Q_{i-1,j}$ and $Q_{i,j}$, let $W^p_{i-1/2,j} = \gamma_p A_{(\sigma_A)_p}$ (or $\delta_p B_{(\sigma_B)_p}$). Then we have
\begin{align}
    \theta_{i-1/2,j}^p = \frac{W^p_{I-1/2,j} \cdot W^p_{i-1/2,j}}{||W^p_{i-1/2,j}||^2}\\
    \tilde{\gamma}_p = \phi(\theta_{i-1/2,j}^p) \gamma_p,
\end{align}
where $I = i-1$ if $(\sigma_A)_p>0$ and $I=i+1$ if $(\sigma_A)_p <0$, and $\phi:\mathbb{R}\rightarrow \mathbb{R}$ is a typical limiting function, such as minmod, MC, etc.
In our problems, we apply the minmod limiter on all the horizontal and vertical grid edges. The minmod limiter is given by:
\begin{align}
    \phi_{mm}(\theta) = \max (0,\min(1,\theta)).
\end{align}
Minmod limiters are symmetric: $\phi_{mm}(1/\theta)=\phi_{mm}(\theta)/\theta$, a desired property especially as our $V$-barrier problem is symmetric.

We note that the wave limiting method described above is suitable only for Cartesian edges, where the neighboring index $I$ is defined. For the barrier cut edge, we will resort to using gradient limiters as we will need to employ limiters suitable for non-Cartesian cell edges. We discuss this further in the next section.

\section{Special numerical method on cut cells: cell merging} \label{sec3}
\cref{fig:cutcells_ex} shows all the cut cells we will consider. 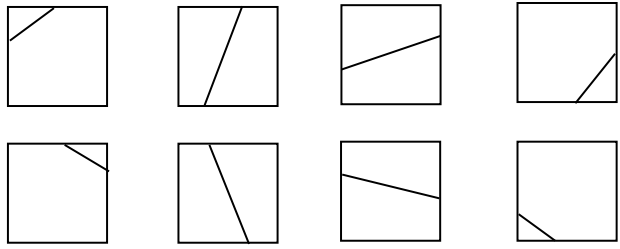
\begin{figure}[h!]
 \tikzset{every picture/.style={line width=0.75pt}} 
\centering
\begin{tikzpicture}[x=0.75pt,y=0.75pt,yscale=-1,xscale=1]

\draw   (338.2,111.1) -- (388.2,111.1) -- (388.2,161.1) -- (338.2,161.1) -- cycle ;
\draw    (171,129) -- (193.2,112.6) ;
\draw   (256,112) -- (306,112) -- (306,162) -- (256,162) -- cycle ;
\draw    (269.2,161.6) -- (288.2,111.6) ;
\draw   (170,112) -- (220,112) -- (220,162) -- (170,162) -- cycle ;
\draw    (338.2,143.6) -- (388.2,126.6) ;
\draw   (427,110) -- (477,110) -- (477,160) -- (427,160) -- cycle ;
\draw    (456.2,160.6) -- (476.2,135.6) ;
\draw   (170,181) -- (220,181) -- (220,231) -- (170,231) -- cycle ;
\draw    (198.6,181.6) -- (221,195) ;
\draw   (256,181) -- (306,181) -- (306,231) -- (256,231) -- cycle ;
\draw    (271.6,181.6) -- (291.6,231.6) ;
\draw   (338,180) -- (388,180) -- (388,230) -- (338,230) -- cycle ;
\draw    (338.6,196.6) -- (387.6,208.6) ;
\draw   (427,180) -- (477,180) -- (477,230) -- (427,230) -- cycle ;
\draw    (427.6,216.6) -- (446,230) ;
\end{tikzpicture}
\caption{All possible types of cut cells for both upper and lower cut cells for a barrier segment (akin to diagram shown in \cite{causon2000calculation}).}
\label{fig:cutcells_ex}
\end{figure}

In this section, we describe in detail the cut cell method we propose for our problem described in \cref{sec1}. First we discuss how the cell merging method works. Then we describe how the fluctuation at the barrier cut edge is computed in both first and second order. Then we discuss how the rest of the cut cell edges are computed in both first and second order.

\subsection{Cell merging}
Cell merging works by absorbing a small cell into a larger neighboring cell, and considering the larger combined cell as its own cell. This way, the CFL limitation on the small cell is avoided, as the combined cell will be large enough to take a full time step with the numerical update.

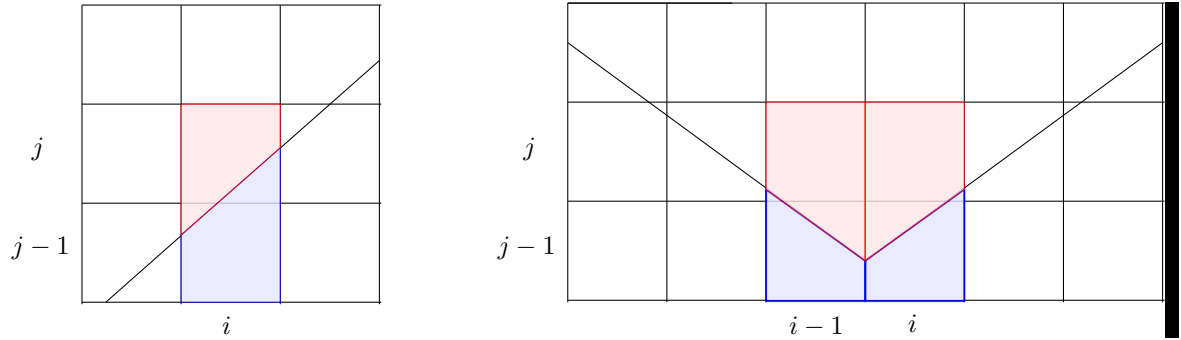
\begin{figure}[h!]
    \centering

\begin{tikzpicture}[x=0.75pt,y=0.75pt,yscale=-1,xscale=1]

\draw  [draw opacity=0] (81.6,81.6) -- (232.6,81.6) -- (232.6,232.6) -- (81.6,232.6) -- cycle ; \draw   (81.6,81.6) -- (81.6,232.6)(131.6,81.6) -- (131.6,232.6)(181.6,81.6) -- (181.6,232.6)(231.6,81.6) -- (231.6,232.6) ; \draw   (81.6,81.6) -- (232.6,81.6)(81.6,131.6) -- (232.6,131.6)(81.6,181.6) -- (232.6,181.6)(81.6,231.6) -- (232.6,231.6) ; \draw    ;
\draw  (231.6,109.6) -- (93.6,231.6); 
\draw [blue,fill=blue!10,opacity=0.8]  (181.6,153.6) -- (181.6,231.6) -- (131.6,231.6) -- (131.6,197.6) -- cycle ;
\draw [red,fill=red!10,opacity=0.8,thin] (131.6,197.6)--(181.6,153.6) -- (181.6,131.6) -- (131.6,131.6) -- (131.6,197.6) -- cycle ;
\draw  [draw opacity=0] (326.6,80.6) -- (627.6,80.6) -- (627.6,231.6) -- (326.6,231.6) -- cycle ; \draw   (326.6,80.6) -- (326.6,231.6)(376.6,80.6) -- (376.6,231.6)(426.6,80.6) -- (426.6,231.6)(476.6,80.6) -- (476.6,231.6)(526.6,80.6) -- (526.6,231.6)(576.6,80.6) -- (576.6,231.6)(626.6,80.6) -- (626.6,231.6) ; \draw   (326.6,80.6) -- (627.6,80.6)(326.6,130.6) -- (627.6,130.6)(326.6,180.6) -- (627.6,180.6)(326.6,230.6) -- (627.6,230.6) ; \draw    ;
\draw    (326.6,100.6) -- (476.6,210.6) ;

\draw    (626.6,100.6) -- (476.6,210.6) ;
\draw [blue, fill=blue!10,opacity=0.8, thick](476.6,231.0)--(476.6,210.6) --(526.6,174.6)-- (526.6,231.0)--(476.6,231.0);
\draw [blue, fill=blue!10,opacity=0.8, thick](476.6,231.0)--(476.6,210.6) --(426.6,174.6)-- (426.6,231.0)--(476.6,231.0);
\draw [red, fill=red!10,opacity=0.8, thin](476.6,210.6) --(426.6,174.6)-- (426.6,130.6)--(476.6,130.6)--(476.6,210.6);
\draw [red, fill=red!10,opacity=0.8, thin](476.6,210.6) --(526.6,174.6)-- (526.6,130.6)--(476.6,130.6)--(476.6,210.6);
\draw (326.6,80.6)--(409.6,80.6);
\draw (55,147.4) node [anchor=north west][inner sep=0.75pt]    {$j$};
\draw (45,197.4) node [anchor=north west][inner sep=0.75pt]    {$j-1$};

\draw (151,237.4) node [anchor=north west][inner sep=0.75pt]    {$i$};
\draw (290,197.4) node [anchor=north west][inner sep=0.75pt]    {$j-1$};
\draw (303,147.4) node [anchor=north west][inner sep=0.75pt]    {$j$};

\draw (497,236.4) node [anchor=north west][inner sep=0.75pt]    {$i$};
\draw (437,237.4) node [anchor=north west][inner sep=0.75pt]    {$i-1$};

\end{tikzpicture}

    \caption{Merged cells using normal neighboring cells for model problems. On the left, we show upper neighboring cell (in red) for cell ($i,j-1$) and lower neighboring cell (in blue) for cell $(i,j)$ and on the right, we show the same for $V$ barrier example for cells ($i-1,j-1$) and ($i,j-1$) and for cells ($i-1,j$) and ($i,j$).}
    \label{fig:merging}
\end{figure}

The first step in the cell merging method is to identify which cells need merging. We set the area threshold to be at $0.5 \Delta x \Delta y$, such that any cell whose area is less than this threshold is a cell that needs to be merged. Second, we need to identify the neighboring cell with which to merge the small cell. In all cases, the neighboring cell will be \emph{the normal} neighboring cell. For our problems, this means the cell directly above the small cell (for small cells above barrier) or directly below the small cell (for small cells below the barrier) as shown in \cref{fig:merging}. The merged cell will then have a volume weighted average of the comprising cells.

The third step is to update the merged cell. We adopt the style of using every edge as done in \cite{2006-Chung-p607} to update our merged cells. This not only gives a conservative method but gives greater accuracy as we do not average the fluctuation over a larger combined edge of a merged cell. We describe how to compute the fluctuation at each edge of the cut cells in the next subsections. Note that if we use a first order fluctuation at each edge, we get a first order cell merging method. If we use a second order fluctuation at each edge, we get a second order cell merging method. That is, cell merging is simply a way to circumvent the CFL restriction and depends on the computational method for fluctuation calculations to either be a low resolution method or a high one.

Final step is to use the updated merged cell to update both the small cell and the merged neighboring cell. If $Q^{n+1}_M$ is the updated merged cell average resulting from merging cells $(i,j)$ and $(i,j+1)$, then we have $Q^{n+1}_{i,j} = Q^{n+1}_{i,j+1} = Q_M^{n+1}$. This is because in essence the cell merging method considers the merged cell as one big cell.

\subsection{Wave redistribution at barrier edge}
\label{WR}
Now we discuss how to set the fluctuation at the barrier edge. We will first define the Riemann problem and the first order fluctuations at the barrier edge, and later in \cref{cut_2nd} we will describe how to do a second order correction at the barrier edge.

\subsubsection{First order}
\begin{figure}[h!]
\centering

\begin{tikzpicture}
\draw (-1.3,-2.1) -- (0.7,-2.1);
\draw (-1.3,-0.1) -- (0.7,-0.1);
\draw (-1.3,-0.1) -- (-1.3,-2.1);
\draw (0.7,-0.1) -- (0.7,-2.1);
\draw[red] (-1.3,-2.1) -- (0.7,-0.1);
\draw[->,very thick,purple] (-0.3,-1.1)--(0.2,-1.6);
\draw[->,very thick,purple] (-0.3,-1.1)--(0.2,-0.6);
\end{tikzpicture}
\caption{Rotation vectors used to rotate states: orthonormal pair.}
\label{fig:rotvec}
\end{figure}
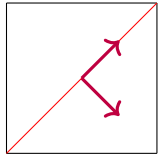

On a cut cell $(i,j)$, we simply use the small cells' state averages as they are on either side of the barrier edge and rotate them as follows:

\begin{align}
    \tilde{Q}^h_{i,j} = R_{i,j} Q^h_{i,j},
\end{align}
where $h=L$ or $U,$ and
\begin{align*}
    R_{i,j} = \begin{bmatrix}
    1 & 0 & 0 \\
    0 & \hat{n}_1 & \hat{n}_2 \\
    0 & \hat{t}_1 & \hat{t}_2
    \end{bmatrix}.
\end{align*}
The rotation vectors are simply the orthonormal pair with respect to the barrier \cref{fig:rotvec}.

We remark that this usage of state variables instead of normal extension works well especially for the V barrier problem because it avoids having to compute extension averages that overlap with each other at the central tip.

Once we have rotated the states, we apply the \emph{wave redistribution} algorithm. The wave redistribution is so called because it is the redistributing of waves that arise from two ghost Riemann problems. At the barrier edge, because we have a zero-width line interface in lieu of actual state representing the barrier, we resort to solving ghost problems. The ghost cell $Q^*$ that is introduced in between the two rotated states is a cell with $B^*=\beta$, the barrier height, and $Q^*=0$ if the heights on both sides are less than the barrier height, to mimic a dry state on top of the barrier. In the case that both water heights of the normal extensions are higher than the barrier height, we choose for our ghost height $H^*$ to be the minimum of the two heights minus $\beta$ and choose the minimum of the two momentum on either side to be our ghost momentum $HU^*,\, HV^*$.

With the ghost cell, we have the following decompositions to perform:
\begin{align}
    f(\tilde{Q}^{L}_{i,j})-f(Q^*) - \Psi_L = A_L\gamma \\
    f(Q^*)-f(\tilde{Q}^{U}_{i,j}) - \Psi_U = A_U\delta,
\end{align}
where $A_L$ is the matrix similar to \cref{matrices} for the Riemann problem between $\tilde{Q}^{L}_{i,j}$ and $Q^*$ with associated eigenvalues and eigenvectors $\{\sigma_L^i, \mathbf{r}_L^i\}$ and $A_U$ is that for the Riemann problem between $\tilde{Q}^{U}_{i,j}$ and $Q^*$ with eigenvalues and eigenvectors $\{\sigma_U^i, \mathbf{r}_U^i\}$. The terms $\Psi_L$ and $\Psi_U$ represent the source term vector: $\Psi_L = [0,g\overline{H}_L(B_L-B^*),0]$ and $\Psi_U=[0,g\overline{H}_U(B^*-B_U),0]$, where $\overline{H}_{L/U}$ is $0.5(H_{L/U}+H^*)$. This simple subtraction of source term from the flux differences is one of the key merits of the wave decomposition method \cite{bale2003wave}.

Then wave redistribution sets a new set of eigenvalues and eigenvectors $\{\omega_i, \boldsymbol{\rho}_i\}$ as
\begin{align}
    &\omega_i = \frac{1}{2} ( \sigma_U^i + \sigma_L^i) \\
    &\boldsymbol{\rho}_i = [1,\omega_i,\frac{1}{2}(\tilde{V}^L + \tilde{V}^U)] \text{ for $i=1,3$} \\
    &\boldsymbol{\rho}_2 = [0,0,1],
\end{align}
where $\tilde{V}^{L/U}$ is the rotated transverse velocity of $\tilde{Q}^{L/U}_{i,j}.$

Finally, wave redistribution solves
\begin{align}
    [A_L | A_U] (\gamma : \delta) = \boldsymbol{\mathrm{P}} \epsilon,
\end{align}
where $[A_L | A_U]$ is the augmented matrix of two matrices $A_L$ and $A_U$, $\gamma : \delta$ is the augmented vector of coefficient vectors $\gamma$ and $\delta$, along the same column axis, and $\boldsymbol{\mathrm{P}}$ is the matrix $[\boldsymbol{\rho}_1, \boldsymbol{\rho}_2, \boldsymbol{\rho}_3]$. The coefficient vector $\epsilon$ is the unknown to be solved for and once solved, we have for our redistributed waves at the barrier edge,
\begin{align}
    A^+\Delta \tilde{Q}_{i,j} = \sum_{p: \omega_p > 0}\epsilon_p \omega_p \boldsymbol{\rho}_p \\
    A^-\Delta \tilde{Q}_{i,j} = \sum_{p: \omega_p< 0}\epsilon_p \omega_p \boldsymbol{\rho}_p.
\end{align}

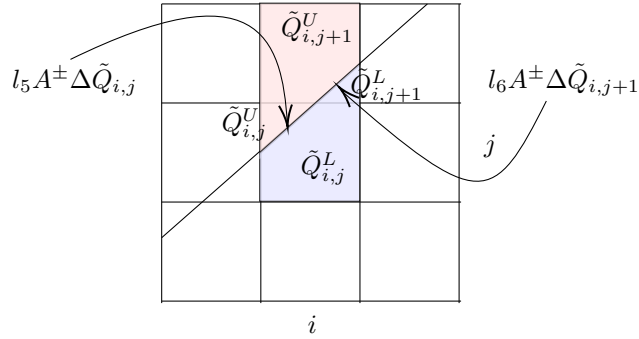
\begin{figure}[h!]
    \centering
    \hfill
\centering

\begin{tikzpicture}[x=0.75pt,y=0.75pt,yscale=-1,xscale=1]

\draw  [draw opacity=0] (231.6,81.6) -- (382.6,81.6) -- (382.6,232.6) -- (231.6,232.6) -- cycle ; \draw   (231.6,81.6) -- (231.6,232.6)(281.6,81.6) -- (281.6,232.6)(331.6,81.6) -- (331.6,232.6)(381.6,81.6) -- (381.6,232.6) ; \draw   (231.6,81.6) -- (382.6,81.6)(231.6,131.6) -- (382.6,131.6)(231.6,181.6) -- (382.6,181.6)(231.6,231.6) -- (382.6,231.6) ; \draw    ;
\draw   (365.6,81.6) -- (231.6,199.6)  ;
\draw [fill=red!10,opacity=0.8] (281,160.4-3.) -- (331.5,115.4-3.) -- (331.5,85-3.8) -- (281,85-3.8) -- (281,163.4-3.8);
\draw [fill=blue!10,opacity=0.8] (281,160.4-3.8) -- (331.5,115.4-3.8) -- (331.5,185-3.8) -- (281,185-3.8) -- (281,163.4-3.8);
\draw    (185.6,109.6) .. controls (278.83,64.06) and (292.35,96.93) .. (294.99,142.22) ;
\draw [shift={(295.1,143.6)}, rotate = 264.40999999999997] [color={rgb, 255:red, 0; green, 0; blue, 0 }  ][line width=0.75]    (10.93,-3.29) .. controls (6.95,-1.4) and (3.31,-0.3) .. (0,0) .. controls (3.31,0.3) and (6.95,1.4) .. (10.93,3.29)   ;
\draw    (425.6,130.6) .. controls (398.83,184.06) and (392.35,186.93) .. (321.99,125.22) ;
\draw [shift={(321.1,123.6)}, rotate = 45] [color={rgb, 255:red, 0; green, 0; blue, 0 }  ][line width=0.75]    (10.93,-3.29) .. controls (6.95,-1.4) and (3.31,-0.3) .. (0,0) .. controls (3.31,0.3) and (6.95,1.4) .. (10.93,3.29)   ;
\draw (393,146.4) node [anchor=north west][inner sep=0.75pt]    {$j$};
\draw (304,237.4) node [anchor=north west][inner sep=0.75pt]    {$i$};
\draw (155,110.4) node [anchor=north west][inner sep=0.75pt]    {$l_{5} A^{\pm} \Delta \tilde{Q}_{i,j}$};
\draw (395,110.4) node [anchor=north west][inner sep=0.75pt]    {$l_{6} A^{\pm} \Delta \tilde{Q}_{i,j+1}$};
\draw (300.6,154) node [anchor=north west][inner sep=0.75pt]    {$\tilde{Q}^{L}_{i,j}$};
\draw (260.6,132) node [anchor=north west][inner sep=0.75pt]    {$\tilde{Q}^{U}_{i,j}$};
\draw (325.6,113) node [anchor=north west][inner sep=0.75pt]    {$\tilde{Q}^{L}_{i,j+1}$};
\draw (290.6,84) node [anchor=north west][inner sep=0.75pt]    {$\tilde{Q}^{U}_{i,j+1}$};
\end{tikzpicture}
\caption{The rotated averages $\tilde{Q}_{i,j}^{U/L}$ and $\tilde{Q}_{i,j+1}^{U/L}$ are used to produce waves at the barrier edges of cut cell $(i,j)$ and cell $(i,j+1)$ to update the merged cells in blue and red. The red merging is for upper cell $(i,j)$ and blue for lower cell $(i,j+1).$}
    \label{fig:barst}
\end{figure}

Once we have computed the positive and negative fluctuations arising from the rotated Riemann problem, we follow the algorithm in \cite{calhoun2008logically} and rotate them back into the original coordinate directions by the following linear transformation:
\begin{align}
A^{\pm} \Delta Q_{i,j} = R_{i,j}^T A^{\pm}\Delta \tilde{Q}_{i,j}.
\end{align}
These waves are then weighted by the length of the barrier cut edge, as shown in \cref{fig:barst}. Note that as shown in \cref{fig:barst}, because we use every piecewise edge of merged cells, we may need to compute fluctuations at two barrier edges to update a single merged cell.

\subsubsection{Second order}\label{cut_2nd}
To apply second order method for the barrier edge fluctuation, we resort to using gradient approximations since this method is conceptually much clearer than using the wave propagation second order method, which works on a logically quadrilateral grid with a clear, linear 3 point stencil (i.e. index $I$). However, for a Riemann problem at a barrier edge, it is not clear what $I=i+1$ or $i-1$ should indicate. Using gradients however, it is clear that we should linearly approximate the solution at the midpoint of the barrier edge and compute the Riemann problem there.

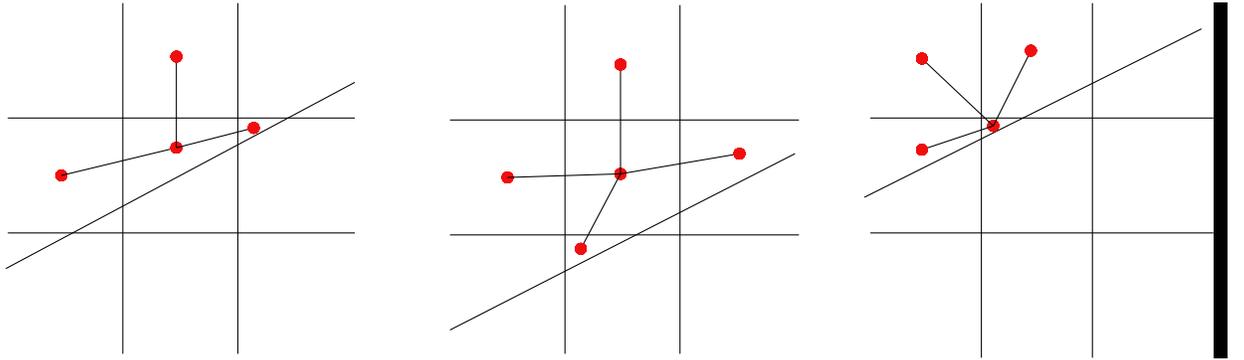
\begin{figure}[h!]
    \centering
\begin{tikzpicture}[x=0.75pt,y=0.75pt,yscale=-1,xscale=1]

\draw  [draw opacity=0]  (186-30,52) -- (186-30,229) -- (11-30,229) -- cycle ; \draw    (11-30,229)(69-30,52) -- (69-30,229)(127-30,52) -- (127-30,229)(185-30,52)  ; \draw    (186-30,52)(11-30,110) -- (186-30,110)(11-30,168) -- (186-30,168)(11-30,226)  ; \draw    ;
\draw    (10-30,186) -- (186-30,92) ;
\draw [color={rgb, 255:red, 241; green, 14; blue, 14 } ] plot[color=red,mark=*,mark size=2pt] (8,139);
\draw (8,139) -- (66,125);
\draw [color={rgb, 255:red, 241; green, 14; blue, 14 } ] plot[color=red,mark=*,mark size=2pt] (66,125);
\draw (66,79) -- (66,125);
\draw (105,115) -- (66,125);
\draw [color={rgb, 255:red, 241; green, 14; blue, 14 } ] plot[color=red,mark=*,mark size=2pt] (66,79);
\draw [color={rgb, 255:red, 241; green, 14; blue, 14 } ] plot[color=red,mark=*,mark size=2pt] (105,115);

\draw  [draw opacity=0]  (410-30,53) -- (410-30,229) -- (234-30,229) -- cycle ; \draw    (234-30,229)(292-30,53) -- (292-30,229)(350-30,53) -- (350-30,229)(408-30,53) ; \draw  (410-30,53)(234-30,111) -- (410-30,111)(234-30,169) -- (410-30,169)(234-30,227)  ; \draw    ;
\draw    (234-30,217) -- (408-30,128) ;
\draw [color={rgb, 255:red, 241; green, 14; blue, 14 } ] plot[color=red,mark=*,mark size=2pt] (290,138.18);
\draw (290,138.18) -- (290,83);
\draw [color={rgb, 255:red, 241; green, 14; blue, 14 } ] plot[color=red,mark=*,mark size=2pt] (233,140);
\draw (233,140) -- (290,138.18);
\draw (350,128) -- (290,138.18);
\draw (270,176) -- (290,138.18);
\draw [color={rgb, 255:red, 241; green, 14; blue, 14 } ] plot[color=red,mark=*,mark size=2pt] (290,83);
\draw [color={rgb, 255:red, 241; green, 14; blue, 14 } ] plot[color=red,mark=*,mark size=2pt] (350,128);
\draw [color={rgb, 255:red, 241; green, 14; blue, 14 } ] plot[color=red,mark=*,mark size=2pt] (270,176);
\draw  [draw opacity=0]  (619-30,52) -- (619-30,231)  -- cycle ; \draw   (446-30,231)(502-30,52) -- (502-30,231)(558-30,52) -- (558-30,231)(614-30,52)  ; \draw    (619-30,52)(446-30,110) -- (619-30,110)(446-30,168) -- (619-30,168)(446-30,220) ; \draw    ;
\draw    (443-30,150) -- (613-30,65) ;
\draw [color={rgb, 255:red, 241; green, 14; blue, 14 } ] plot[color=red,mark=*,mark size=2pt] (478,114);
\draw (442,80)--(478,114);
\draw (478,114)--(442,126);
\draw (497,76)--(478,114);
\draw [color={rgb, 255:red, 241; green, 14; blue, 14 } ] plot[color=red,mark=*,mark size=2pt] (497,76);
\draw [color={rgb, 255:red, 241; green, 14; blue, 14 } ] plot[color=red,mark=*,mark size=2pt] (442,80);
\draw [color={rgb, 255:red, 241; green, 14; blue, 14 } ] plot[color=red,mark=*,mark size=2pt] (442,126);

\end{tikzpicture}
    \caption{Stencils used for approximating gradients on different types of cut cells. The cut cell of interest is the upper cut cells at the center.}
    \label{fig:stencil_grad}
\end{figure}

To approximate the gradient, we use the following stencils depending on which cut cell we are approximating the gradient for. There are three types of stencils we use, as shown in \cref{fig:stencil_grad}. In the diagram we focus on the upper cut cells at the center. However, the same principle applies to the lower cut cells and also the cut cells with barrier of slope $m>1.$ We use the regular five point stencil with the exception of cells that the barrier blocks the central cell from accessing.

To limit the gradients we use the commonly used Barth-Jespersen limiter. This limiter has the advantage of keeping the gradients low such that the reconstructed value will always be in between the maximum and minimum values over the stencil. This assures nice properties such as positivity for the height variable, given everywhere else height is positive. Also, it avoids oscillations by not introducing new maxima.

The Barth-Jespersen limiter $\alpha$ is computed by taking the minimum of $\alpha_{N_{k}}$, where $N_k$ is a member of stencil $N_{i,j}=\{N_k\}$:
\begin{align}
    \alpha = \min\{ \alpha_{N_k}\},
\end{align}
where
\begin{align}
    \alpha_{N_k} = \begin{cases} \min (1, (M_{i,j} - Q_{i,j}^n)/(Q_{N_k}^n - Q_{i,j}^n) ) &\mbox{if } Q^n_{N_k} - Q^n_{i,j} > 0 \\
    \min (1,(m_{i,j} - Q_{i,j}^n)/(Q_{N_k}^n - Q_{i,j}^n) ) &\mbox{if } Q^n_{N_k} - Q^n_{i,j} < 0\\
    1 &\mbox{if } Q^n_{N_k} - Q_{i,j}^n = 0, \end{cases}
\end{align}
where $M_{i,j}$ denotes the maximum solution value over stencil and $m_{i,j}$ the minimum.
We simply multiply the limiter to $\nabla Q^n_{i,j}$ to give $\nabla \tilde{Q}^n_{i,j} = \alpha \nabla Q^n_{i,j}$. Note that since our $Q \in \mathbb{R}^3$, we apply the limiter to each variable in $Q$.

With the limited gradients calculated on either side of the cut, we compute the linearly reconstructed value $Q$ at the barrier edge from both sides and apply wave redistribution on those values as described in \cref{WR}.

\subsection{At non-barrier edges}
At every other edge, which are the vertical and horizontal edges of the cut cell, we employ the base method as described in \cref{wp_2nd} with two caveats.
\begin{figure}[h!]
    \centering
\begin{tikzpicture}[x=0.75pt,y=0.75pt,yscale=-1,xscale=1]

\draw  [draw opacity=0] (241.6,74.6) -- (392.6,74.6) -- (392.6,225.6) -- (241.6,225.6) -- cycle ; \draw   (241.6,74.6) -- (241.6,225.6)(291.6,74.6) -- (291.6,225.6)(341.6,74.6) -- (341.6,225.6)(391.6,74.6) -- (391.6,225.6) ; \draw   (241.6,74.6) -- (392.6,74.6)(241.6,124.6) -- (392.6,124.6)(241.6,174.6) -- (392.6,174.6)(241.6,224.6) -- (392.6,224.6) ; \draw    ;
\draw    (379.6,74.6) -- (241.6,196.6) ;
\draw [fill=red!10,opacity=0.8] (291,153.4) -- (341.5,108.4) -- (341.5,75) -- (291,75) -- (291,153.4);
\draw [fill=blue!10,opacity=0.8] (291,153.4) -- (341.5,108.4) -- (341.5,175) -- (291,175) -- (291,153.4);
\draw    (264.6,244.4) .. controls (246.78,202.82) and (247.58,193.58) .. (290.68,162.94) ;
\draw [shift={(292,162)}, rotate = 504.73] [color={rgb, 255:red, 0; green, 0; blue, 0 }  ][line width=0.75]    (10.93,-3.29) .. controls (6.95,-1.4) and (3.31,-0.3) .. (0,0) .. controls (3.31,0.3) and (6.95,1.4) .. (10.93,3.29)   ;
\draw    (407.6,170.4) .. controls (426.5,102.74) and (391.95,81.61) .. (353.88,107.37) ;
\draw [shift={(353,109)}, rotate = 323.98] [color={rgb, 255:red, 0; green, 0; blue, 0 }  ][line width=0.75]    (10.93,-3.29) .. controls (6.95,-1.4) and (3.31,-0.3) .. (0,0) .. controls (3.31,0.3) and (6.95,1.4) .. (10.93,3.29)   ;
\draw    (315.6,46.4) .. controls (255.82,-2.6) and (281.52,79.97) .. (320.57,80) ;
\draw [shift={(322.6,81.8)}, rotate = 200.31] [color={rgb, 255:red, 0; green, 0; blue, 0 }  ][line width=0.75]    (10.93,-3.29) .. controls (6.95,-1.4) and (3.31,-0.3) .. (0,0) .. controls (3.31,0.3) and (6.95,1.4) .. (10.93,3.29)   ;
\draw    (218.6,112.4) .. controls (219.59,122.3) and (222.15,169.84) .. (272.07,143.23) ;
\draw [shift={(273.6,142.4)}, rotate = 511.08] [color={rgb, 255:red, 0; green, 0; blue, 0 }  ][line width=0.75]    (10.93,-3.29) .. controls (6.95,-1.4) and (3.31,-0.3) .. (0,0) .. controls (3.31,0.3) and (6.95,1.4) .. (10.93,3.29)   ;

\draw (396,140.4) node [anchor=north west][inner sep=0.75pt]    {$j$};
\draw (310,226.4) node [anchor=north west][inner sep=0.75pt]    {$i$};
\draw (295,152.4) node [anchor=north west][inner sep=0.75pt]    {$l_{1}$};
\draw (343,107.4) node [anchor=north west][inner sep=0.75pt]    {$l_{2}$};
\draw (277,128.4) node [anchor=north west][inner sep=0.75pt]    {$l_{3}$};
\draw (328,80.8) node [anchor=north west][inner sep=0.75pt]    {$l_{4}$};
\draw (230,254.4) node [anchor=north west][inner sep=0.75pt]    {$l_{1} A^{+} \Delta Q_{i-1/2,j}^{L}$};
\draw (375,173.4) node [anchor=north west][inner sep=0.75pt]    {$l_{2} A^{-} \Delta Q_{i+1/2,j+1}^{L}$};
\draw (296,50.4) node [anchor=north west][inner sep=0.75pt]    {$l_{4} A^{-} \Delta Q_{i+1/2,j+1}^{U}$};
\draw (183,92.4) node [anchor=north west][inner sep=0.75pt]    {$l_{3} A^{+} \Delta Q_{i-1/2,j}^{U}$};

\end{tikzpicture}
    \caption{The length ($l_i$) weighted waves on cut edges used to update merged cell for upper cut cell $(i,j)$ in red and for lower cut cell $(i,j+1)$ in blue. The fluctuations are also computed as usual at the uncut edges (e.g. bottom edge and right edge of cell $(i,j)$ for the blue merged cell). }
    \label{fig:dblstate}
\end{figure}
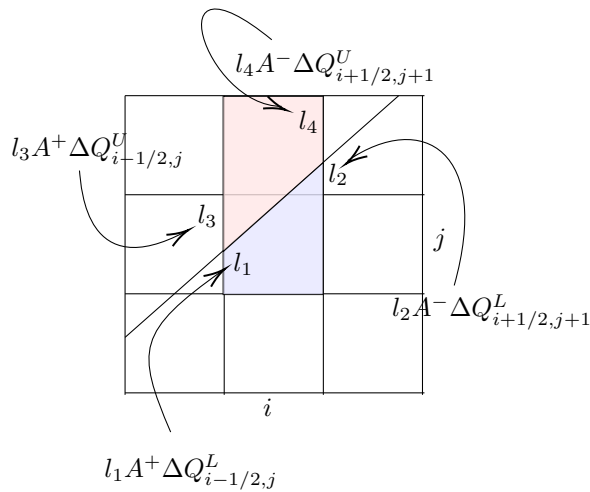
\subsubsection{First order}

In the first order, the $A^{\pm}\Delta Q_{i-1/2,j}$ and $B^{\pm}\Delta Q_{i,j-1/2}$ are computed exactly in the same way as described in \cref{wp_1st}.
Because of the shorter length of the cut cell edges, however, we need to \emph{weight} the fluctuations that arise by the length of the cut cell edge to $\Delta x$ \cite{CALHOUN2000143}, just as was done with the fluctuation at the barrier edges. This is the first caveat and shown diagrammatically in \cref{fig:dblstate}.

All in all, the update formula for the lower cut cell in \cref{fig:dblstate} in first order will look as follows:
\begin{align}
    Q^{L,n+1}_{i,j} & = \frac{(\alpha^L_{i,j} Q^{L,n}_{i,j} + \alpha^L_{i,j+1} Q^{L,n}_{i,j+1})}{\alpha^L_{i,j}+\alpha^L_{i,j+1}} - \frac{\Delta t}{\alpha^L_{i,j}+\alpha^L_{i,j+1}}(l_{i,j}A^+ \Delta Q_{i,j} + l_{i,j+1} A^+\Delta Q_{i,j+1} \nonumber \\
   & + l^L_{i-1/2,j}A^+\Delta Q^L_{i-1/2,j} + l^L_{i+1/2,j}A^-\Delta Q^L_{i+1/2,j} \nonumber \\& + l^L_{i+1/2,j+1} A^-\Delta Q^L_{i+1/2,j+1} + l^L_{i,j-1/2}B^+\Delta Q^L_{i,j-1/2}),
\end{align}
where $\alpha^{U/L}_{i,j}, \, l_{i,j}$ denote the area of cut cell and the length of the barrier edge, and $l^{U/L}_{i\pm 1/2,j} $ and $l^{U/L}_{i,j\pm 1/2}$ represent the lengths of the vertical and horizontal edges of the cut cell, respectively. The upper cut cell is updated in a similar manner.

\subsubsection{Second order}
\begin{figure}
    \centering[h!]
\begin{tikzpicture}[x=0.75pt,y=0.75pt,yscale=-1,xscale=1]

\draw  [draw opacity=0] (100,70) -- (257,70) -- (257,225) -- (100,225) -- cycle ; \draw   (154,70) -- (154,225)(208,70) -- (208,225) ; \draw   (100,124) -- (257,124)(100,178) -- (257,178) ; \draw    ;
\draw    (241,60) -- (94,210) ;
\draw [color={rgb, 255:red, 241; green, 14; blue, 14 }] plot[color=red,mark=x,mark size=4pt] (154,135);
\draw [color={rgb, 255:red, 14; green, 27; blue, 241 }] plot[color=red,mark=x,mark size=4pt] (180,179);
\draw [color={rgb, 255:red, 241; green, 14; blue, 14}] plot[color=red,mark=x,mark size=4pt] (165,124);
\draw [color={rgb, 255:red, 14; green, 27; blue, 241}] plot[color=red,mark=x,mark size=4pt] (208,150);
\draw  [draw opacity=0] (404,66) -- (554,66) -- (554,215) -- (404,215) -- cycle ; \draw   (404,66) -- (404,215)(551,66) -- (551,215) ; \draw   (404,66) -- (554,66)(404,213) -- (554,213) ; \draw    ;
\draw    (516,47) -- (383,173) ;
\draw [color={rgb, 255:red, 241; green, 14; blue, 14}]  (431,106.5) -- (444,101) -- (444,103.75) -- (461.5,103.75) -- (461.5,101) -- (474.5,106.5) -- (461.5,112) -- (461.5,109.25) -- (444,109.25) -- (444,112) -- cycle ;
\draw [color={rgb, 255:red, 241; green, 14; blue, 14}]  (447,95.06) -- (453,84) -- (459,95.06) -- (456,95.06) -- (456,117.19) -- (459,117.19) -- (453,128.25) -- (447,117.19) -- (450,117.19) -- (450,95.06) -- cycle ;
\draw [color={rgb, 255:red, 241; green, 14; blue, 14}]  (378,106.5) -- (391,101) -- (391,103.75) -- (408.5,103.75) -- (408.5,101) -- (421.5,106.5) -- (408.5,112) -- (408.5,109.25) -- (391,109.25) -- (391,112) -- cycle ;
\draw [color={rgb, 255:red, 241; green, 14; blue, 14}]  (447,45.81) -- (453,34.75) -- (459,45.81) -- (456,45.81) -- (456,67.94) -- (459,67.94) -- (453,79) -- (447,67.94) -- (450,67.94) -- (450,45.81) -- cycle ;
\draw (270,147) node [anchor=north west][inner sep=0.75pt]    {$j$};
\draw (177,40.4) node [anchor=north west][inner sep=0.75pt]    {$i$};

\end{tikzpicture}

    \caption{Barrier is blocking access from the upper cut cell $(i,j)$ to its right and bottom neighboring cell, when speed and wave information from Riemann problems there (crossed in blue) are required to perform limiting (\cref{wave_lim}). In such cases, we use the wave and speed that arise from the barrier (denoted by the overlapping vertical and horizontal double-ended arrows in the right figure) to correct the fluctuations at the left and top edge. The wave and speed from the barrier is first computed normal to the barrier then rotated back to $x$ and $y$ direction.}
    \label{fig:lim_corr}
\end{figure}
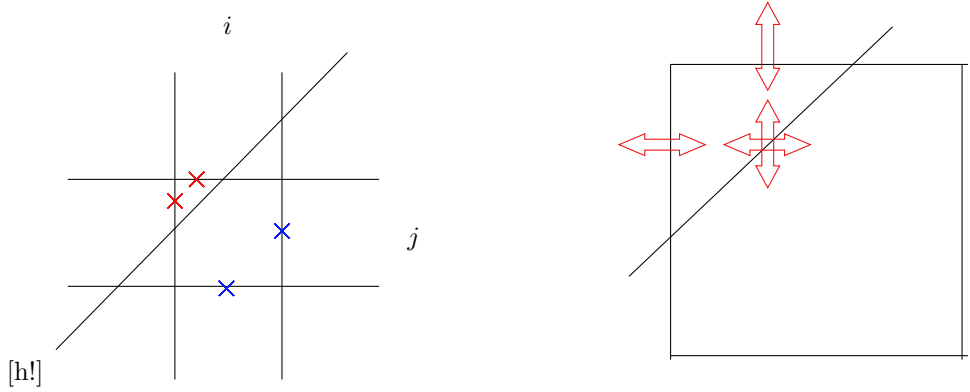

In second order, waves and limited wave corrections are also weighted. However, the second caveat is in computing second order correction in the non-barrier cut cell edges. This is because some cut cells do not have a right or left edge, due to the barrier completely blocking it, when information from the right or left edge is necessary to perform limiting (the index $I$ in \cref{wave_lim}). This is resolved by using the wave and speed arising from the barrier (\cref{fig:lim_corr}).

\section{Numerical examples}  \label{sec4}
All the numerical examples presented below have steady water height of $h=1.2$ and a dam jump of $\Delta h = 0.8$, giving the overall height of the dam break to be $h=2.0$. The barrier height $\beta$ is chosen to be $\beta=1.5$. These are chosen to test overtopping of wave over the barrier. For the boundary conditions, we have one extrapolation boundary condition on the side that the overtopping waves travel, such that they can exit the domain after overtopping the barrier. Furthermore, these computations are done using the second order method, with the Barth Jespersen limiter for gradient at the cut edge and minmod limiter at Cartesian edges.

For comparison, we run simulations with the same initial condition with a mapped grid suited for each barrier, which we go in more detail in the following sections. We observe the 2D color contour plots of water height at specific times from both the CM cut cell method and mapped example runs. Also we observe gauge data, which are time profiles of water height at a specified location in the domain, marked by an asterisk with a number in the figures below.

\subsection{The $20^\circ$ angled barrier}
Here we present the first numerical example with a positively sloped 20$^\circ$ barrier. We place two gauges, one at $(0.5,0.39)$ as before and another at $(0.5,0.8)$ to both capture the reflection and also the overtopped wave further away from the barrier.

\subsubsection{Case of Overtopping}
\begin{figure}[h!]
    \centering
    \includegraphics[scale=0.33]{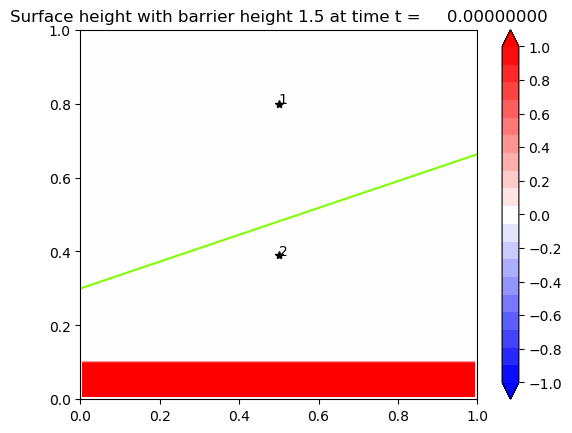}
    \caption{Initial condition for overtopping case: $\beta = 1.5$. Dam height is $2.0$. Grid is $900\times 900.$}
    \label{fig:sot0}
\end{figure}

We show the initial condition in \cref{fig:sot0}.

\begin{figure}[h!]
\begin{subfigure}[b]{0.5\textwidth}
\centering
    \includegraphics[scale=0.33]{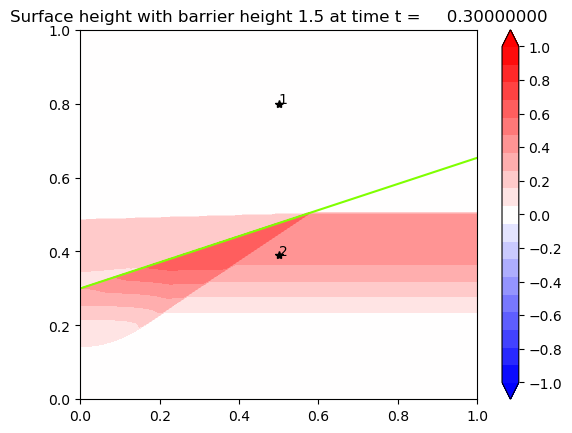}
    \caption{$t=0.3$}
    \label{fig:sot2}
\end{subfigure}
\begin{subfigure}[b]{0.5\textwidth}
\centering
    \includegraphics[scale=0.33]{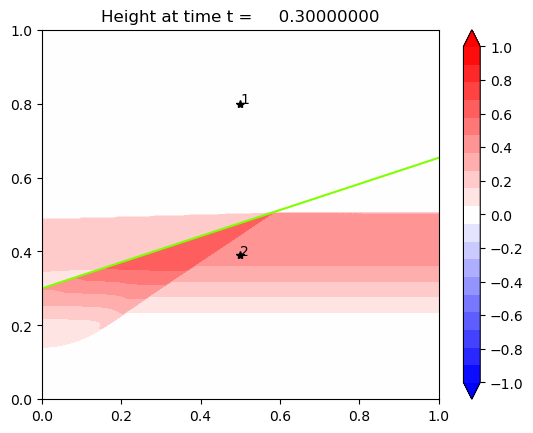}
    \caption{$t=0.3$}
    \label{fig:sot2gc}
\end{subfigure}
\caption{Linear barrier example: CM on left and mapped grid on right.}
        \label{fig:sot2_fig}
\end{figure}

We observe from both CM results and mapped grid results that the wave is abated from the barrier and proceeds in the same direction after overtopping (\cref{fig:sot2_fig}). The reflection wave on the lower side of the barrier is gliding up the linear barrier, while reflecting back in normal direction to the barrier at the same time.

At $t=0.7$, the overtopped wave almost exits the domain at the upper boundary. On the lower side of the barrier, the reflected waves are bouncing around the wall boundary conditions, having reflected from the bottom boundary (the lower left radial wave upward) and from the barrier corner (the upper right radial wave downward).
\begin{figure}[h!]
\begin{subfigure}[b]{0.5\textwidth}

    \centering
    \includegraphics[scale=0.33]{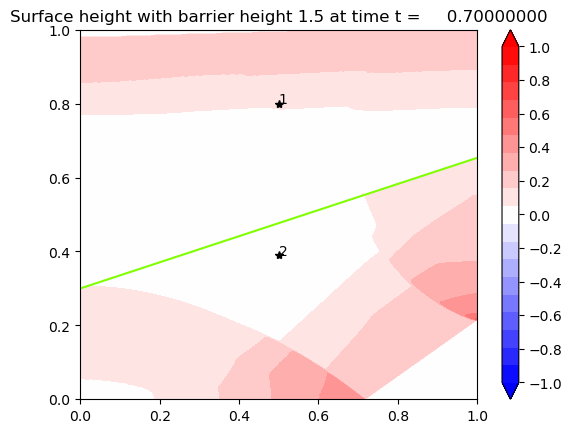}
    \caption{$t=0.7$}
    \label{fig:sot7}
    \end{subfigure}
    \begin{subfigure}[b]{0.5\textwidth}
    \centering
    \includegraphics[scale=0.33]{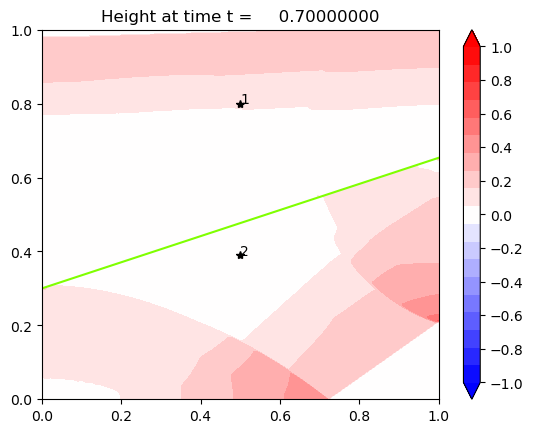}
    \caption{$t=0.7$}
    \label{fig:sot7_gc}
    \end{subfigure}
    \caption{CM (left) and mapped grid (right).}

        \label{fig:sot7_fig}

\end{figure}

An interesting observation to be made is at $t=1.4,$ the doubly-reflected wave (from barrier and bottom edge) glide up the barrier and at the top right corner `pinches up' to achieve overtopping momentum, whereas it does not overtop on the lower part of the barrier.
\begin{figure}[h!]
\begin{subfigure}[b]{0.5\textwidth}

    \centering
    \includegraphics[scale=0.33]{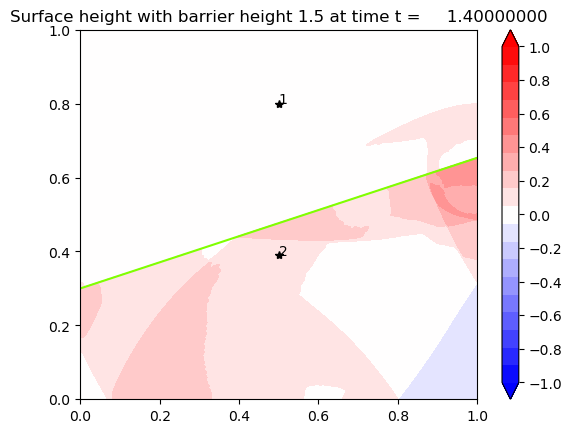}
    \caption{$t=1.4$}
    \label{fig:sot14}
    \end{subfigure}
    \begin{subfigure}[b]{0.5\textwidth}
    \centering
    \includegraphics[scale=0.33]{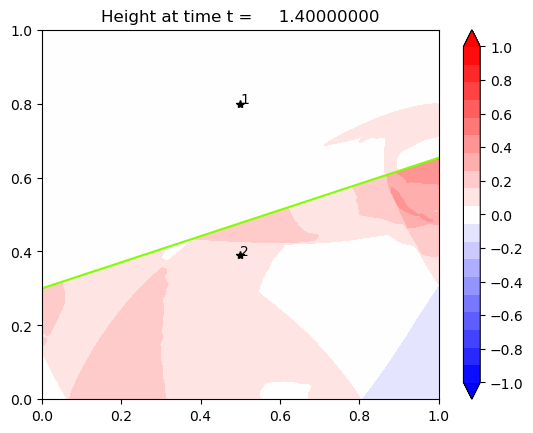}
    \caption{$t=1.4$}
    \label{fig:sot14_gc}
    \end{subfigure}
    \caption{CM (left) and mapped grid (right).}

        \label{fig:sot14_fig}

\end{figure}

\subsubsection{Comparison to mapped grid}
\begin{figure}[h!]
    \centering
    \includegraphics[scale=0.5]{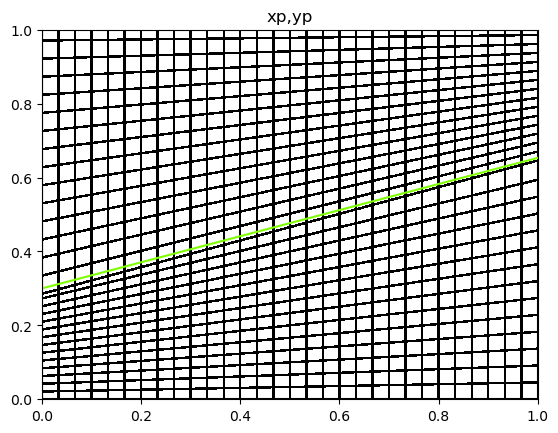}
    \caption{Mapped grid for the $20^\circ$ barrier. Coarsened to $30 \times 30$ to highlight mapping. The lime green line represents the edge where the zero width barrier is located. Along this edge we do wave redistribution in the computationally $y$-direction.}
    \label{fig:mapgridL}
\end{figure}

The way to validate the CM results is by comparing them against a mapped grid example, as shown on the right column in (\cref{fig:sot2_fig}-\cref{fig:sot14_fig}). The mapped grid is shown in \cref{fig:mapgridL}. We transform the computational uniform grid $(x,y)$ into a skewed grid $f_L(x,y)$, with the following mapping $f_L$:
\begin{align}
    f_L(x,y) & = (x, \mu_L (y)), \\
    \mu_L(y) & =\begin{cases} \frac{L(x)}{y^*} y  &\mbox{if } y \in [0,y^*] \\
\frac{1-L(x)}{1-y^*}(y-1) +1 & \mbox{if } y \in [y^*,1]
 \end{cases},
\end{align}
where $L(x)$ is the barrier line equation and $y^*$ is the computational $y$-edge that is mapped to the barrier edge (lime green in \cref{fig:mapgridL}). The computational $y$-edge is taken to be the midpoint of the $y$ coordinates of coordinate 1 and 2 (\cref{fig:2d_setup1}).

This is a linear and affine transformation that stretches and squeezes the upper and lower halves of the computational domain into the physical grid as shown in \cref{fig:mapgridL}. We implement wave redistribution at the computational barrier edge $(y=y^*)$. The finite volume method for mapped grids is explained in \cite{leveque2002finite} and is very much similar to the rotational part of the wave redistribution method explained in \cref{WR}.

\begin{figure}[h!]
\begin{subfigure}[b]{0.5\textwidth}
\centering
    \includegraphics[scale=0.35]{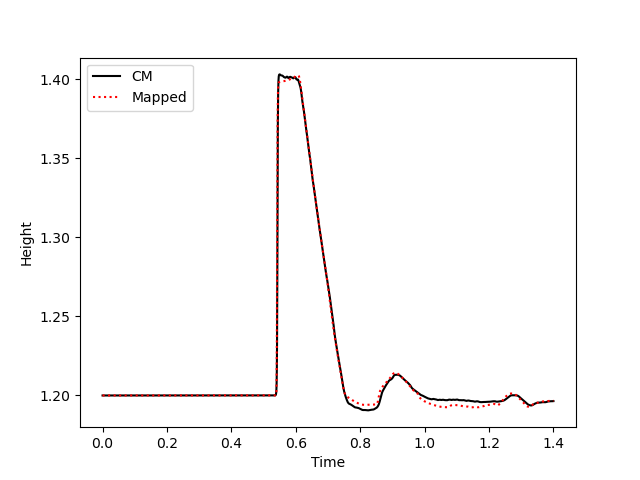}
    \caption{Time profile of Gauge 1 $(0.5,0.8)$.}
    \label{fig:otmp}
\end{subfigure}
\begin{subfigure}[b]{0.5\textwidth}
\centering
    \includegraphics[scale=0.35]{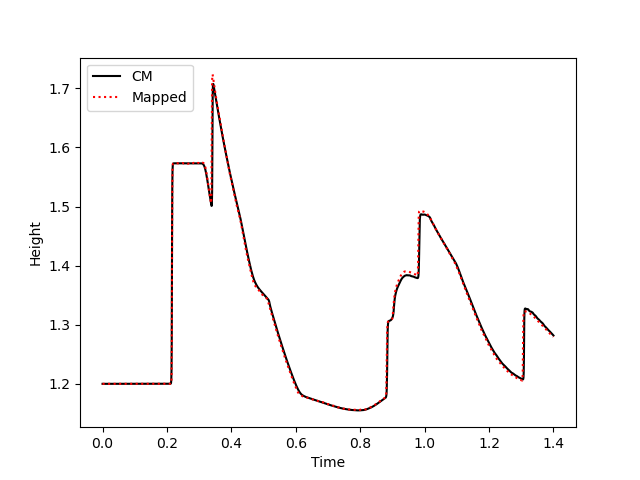}
    \caption{Time profile of Gauge 2 $(0.5,0.39)$.}
    \label{fig:otmp2}
\end{subfigure}
\caption{Gauge comparisons between CM and mapped grid. Results from $900 \times 900$ for CM and from $900 \times 900$ for mapped grid.}
    \label{fig:sotTP}
\end{figure}

We observe the gauge results between the CM and the mapped examples in \cref{fig:sotTP} and see a very good comparison. The gauge results of the mapped example are from a run on $900 \times 900$ grid, and the results of the CM example are also from a $900 \times 900$ grid.
\subsubsection{Convergence}
For convergence studies we observe the convergence of the wave height profile at each gauge point. As our studies are for developing models for protective strategies against storms, we focus not on the numerical results at a slice in the spatial domain at a fixed time, but rather on a slice in the temporal domain at a fixed point in the physical grid. This is often done in storm simulations, to test the storm models' accuracy against real results collected at a specified gauge point off a coast (e.g. Battery Park, NYC).
\begin{table}[h!]
\centering
\begin{tabular}{|c|c|c|c|c|c|}
\hline
~ & ~ & \multicolumn{2}{c|}{$L_1$ Error (1st)} & \multicolumn{2}{c|}{$L_1$ Error (2nd)}\\
$\Delta x$ & $N_x, N_y$ & Gauge 1 & Gauge 2 &Gauge 1 & Gauge 2 \\ \hline
4.e-2   & 25  & 1.03e-2 & 2.45e-2 & 8.95e-3 & 1.37e-2 \\ \hline
2.e-2   & 50  & 3.21e-3 (3.18)& 1.13e-2(2.16) & 2.57e-3 (3.20)& 4.28e-3 (3.48) \\ \hline
1.e-2   & 100 & 9.20e-4 (3.49) & 4.02e-3 (2.81)& 7.27e-4 (4.00)& 1.07e-3 (3.54)\\ \hline
0.666e-2 & 150& 4.13e-4 (2.23) & 2.22e-3 (1.80)& 3.36e-4 (1.49)& 7.17e-4 (2.16)\\ \hline
0.5e-2 & 200 & 2.56e-4 (1.61) & 1.25e-3(1.77)& 1.62e-4 (1.84) & 3.88e-4 (2.08)\\ \hline

\end{tabular}
\caption{$L_1$ errors at Gauge 1 (0.5,0.8) and Gauge 2 (0.5,0.39).}
\label{tab:Lbar_g2}
\end{table}

Also we only study the convergence of the overtopped examples, for similar reasons as abovementioned, namely, that in realistic scenarios the barriers will be overtopped by incoming waves. Furthermore the reflection-only results look promising even in lower resolutions with the matching gauge results and time snapshots of the 2D results.

For our standard solution to compare against, we use the mapped grid example, as it is numerically the same problem as our model problem, and take the results from $900 \times 900$ grid. Then we take the heights at time intervals $\{0.0, \, 0.1, \, ... \, 1.4\}$ with $\Delta t = 0.1$ and compare the CM and mapped results.

\begin{figure}[h!]
    \centering
    \includegraphics[scale=0.5]{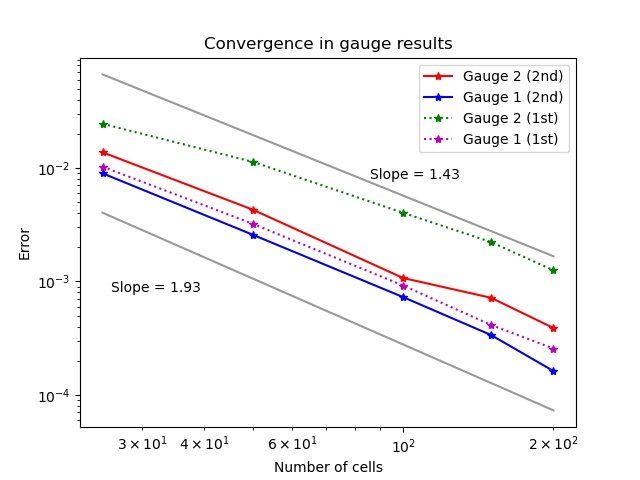}
    \caption{Convergence of gauge profiles.}
    \label{fig:Lconv}
\end{figure}

Shown in \cref{tab:Lbar_g2} are the $L_1$ normed errors for the overtopped wave profile.  We also have the convergence plot in \cref{fig:Lconv}. We observe that the gauge 2 seems to converge little more slowly, especially in the first order method. We attribute this to the fact that it is nearer to the barrier, where the cut cells are. However, in the second order method, gauge 2 performs just as well as gauge 1, because of the gradient reconstruction at the edges and limiting.

\subsection{$117^\circ$ angled V-barrier}
\begin{figure}[h!]
    \centering
    \includegraphics[scale=0.3]{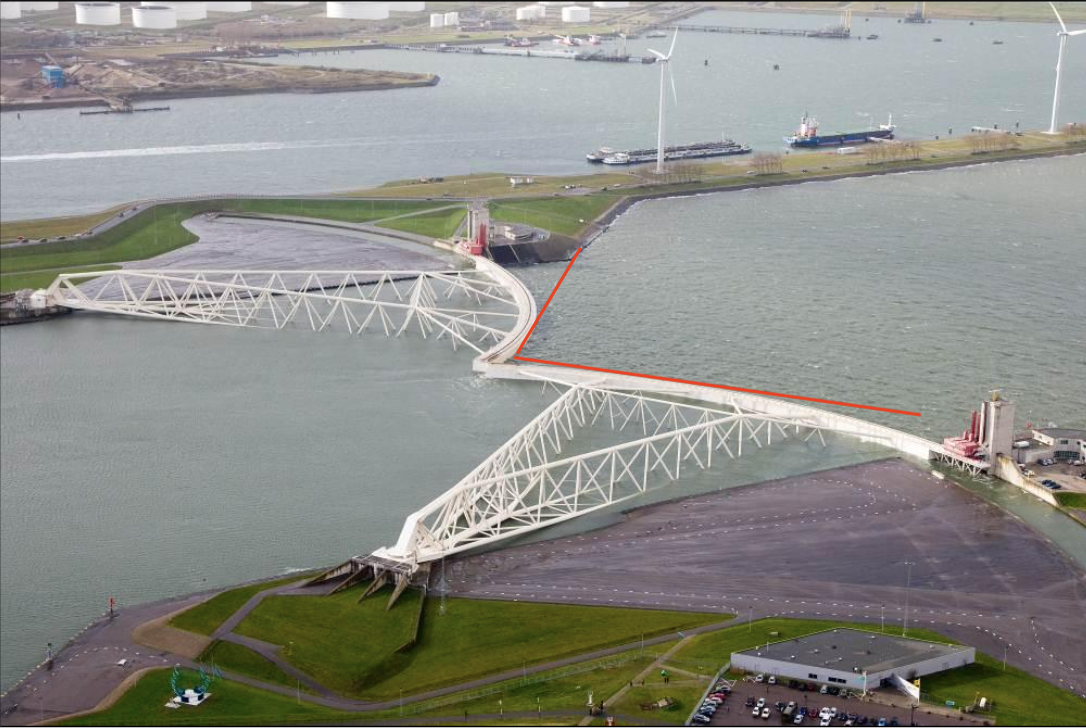}
    \caption{The Maeslant Barrier in the Netherlands. The barrier can approximately be represented as a V-shape (highlighted in red).}
    \label{fig:maeslant}
\end{figure}
\begin{figure}[h!]
    \centering
    \includegraphics[scale=0.33]{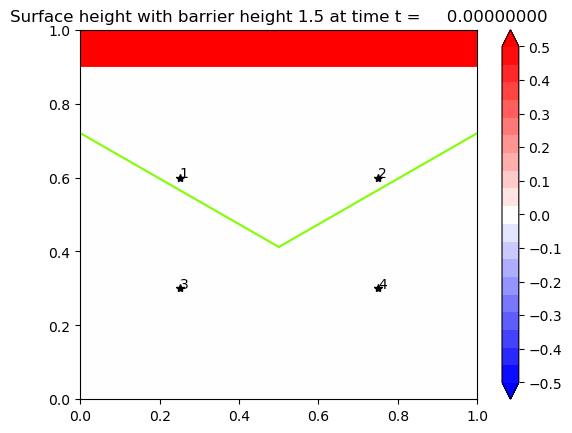}
    \caption{Initial condition for overtopping case: $\beta = 1.5$. Dam height is $2.0$. Grid is $300 \times 300$ for both plots.}
    \label{fig:ot0}
\end{figure}

The inspiration for the V-barrier problem comes from the Maeslant Barrier in the Netherlands (\cref{fig:maeslant}), which closes to form a curved V-shape. We also note that the concave part of the Maeslant barrier is where the incoming wave is expected to hit (i.e. V faces towards the sea). Thus we model our problem such that a dam break initiates a wave from the top side of the V-barrier (\cref{fig:ot0}). Also, in \cite{berger2012simplified}, a similar example is solved for the Euler equations, but with no flux allowed at the V boundary and with computations only for shock wave reflection.

\subsubsection{Case of Overtopping}
\begin{figure}[h!]
\begin{subfigure}[b]{0.5\textwidth}
\centering
    \includegraphics[scale=0.35]{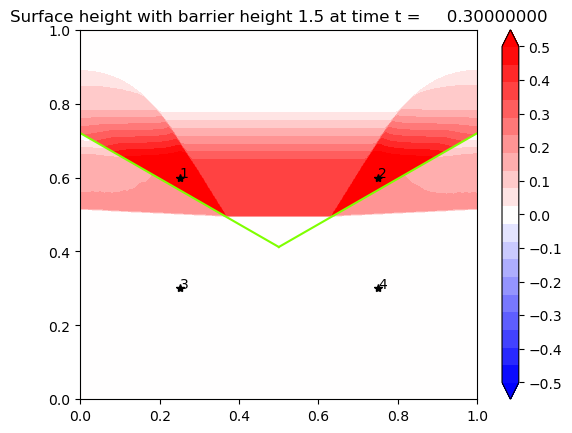}
    \caption{CM: $t=0.3$}
    \label{fig:ot2}
\end{subfigure}
\begin{subfigure}[b]{0.5\textwidth}
\centering
    \includegraphics[scale=0.35]{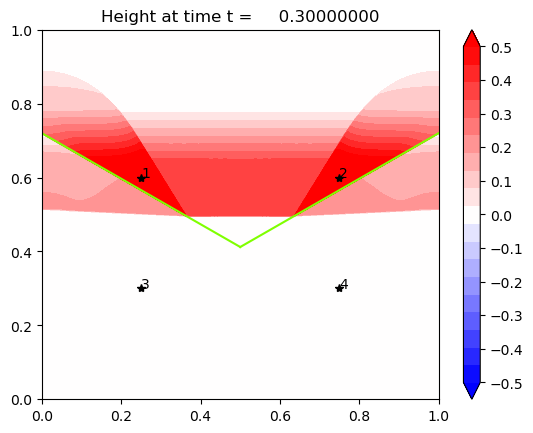}
    \caption{mapped grid: $t=0.3$}
    \label{fig:ot2gc}
\end{subfigure}
\caption{CM Comparison with mapped grid. Both $300 \times 300$ grid. Note the structure of the just overtopped wave.}
        \label{fig:ot2_fig}

\end{figure}

\begin{figure}[h!]
\begin{subfigure}[b]{0.5\textwidth}

    \centering
    \includegraphics[scale=0.35]{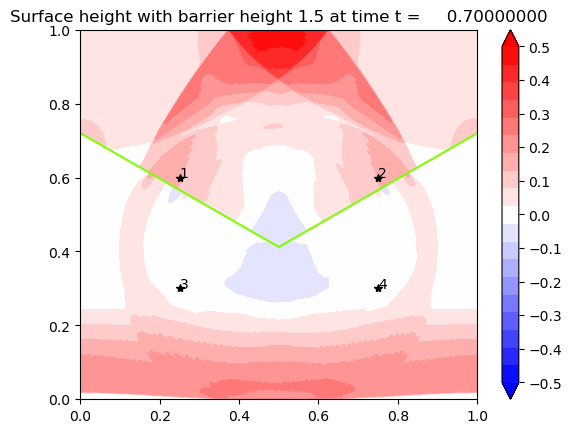}
    \caption{CM: $t=0.7$}
    \label{fig:ot7}
    \end{subfigure}
    \begin{subfigure}[b]{0.5\textwidth}
    \centering
    \includegraphics[scale=0.35]{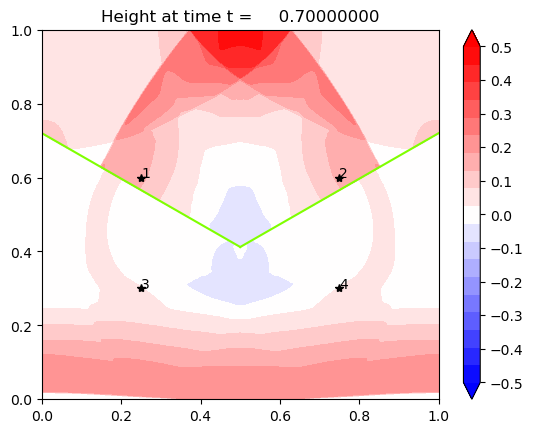}
    \caption{mapped grid: $t=0.7$}
    \label{fig:ot7_gc}
    \end{subfigure}
    \caption{CM Comparison with mapped grid. Note the radially outward moving wave from the center of the V-barrier.}

        \label{fig:ot7_fig}

\end{figure}

\begin{figure}[h!]
\begin{subfigure}[b]{0.5\textwidth}

    \centering
    \includegraphics[scale=0.35]{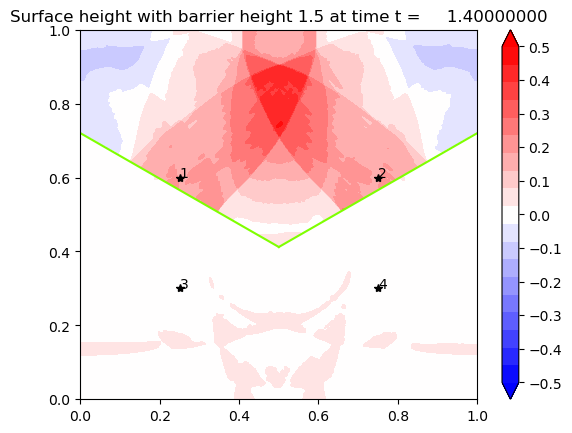}
    \caption{CM: $t=1.4$}
    \label{fig:ot14}
    \end{subfigure}
    \begin{subfigure}[b]{0.5\textwidth}
    \centering
    \includegraphics[scale=0.35]{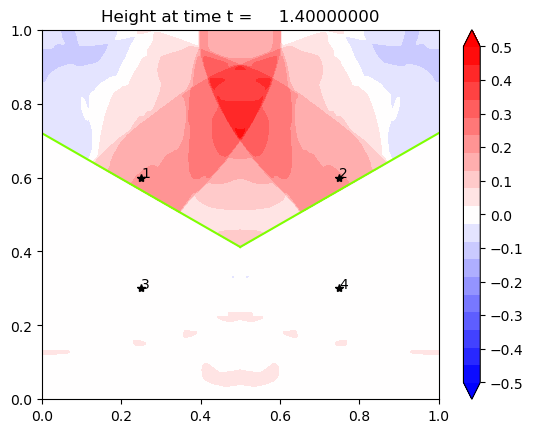}
    \caption{mapped grid: $t=1.4$}
    \label{fig:ot14_gc}
    \end{subfigure}
    \caption{CM Comparison with mapped grid. Note the ``island" of peak at the bottom center.}

        \label{fig:ot14_fig}

\end{figure}

Again we have water height as $1.2$ and the dam jump to be $0.8$ to achieve total dam height of $2.0$. As we shall see, this is enough height to overtop the barrier of height 1.5. We do a comparison against mapped grid results to show accuracy of our CM results. From \cref{fig:ot2_fig} to \cref{fig:ot14_fig}, we can see the similarities of the two 2D plots. At time $t=0.3$, we see the overtopping wave's ``wing"-like structure just below the V-barrier where the amplitude is highest. In \cref{fig:ot7_fig}, we see the overtopping wave moving radially outward from the center, shown both in the CM and mapped grid results. Finally in \cref{fig:ot14_fig} we see the small islands of wave amplitude at the bottom center of the plot in both results.

\begin{figure}[h!]
\begin{subfigure}[b]{0.5\textwidth}
\centering
    \includegraphics[scale=0.35]{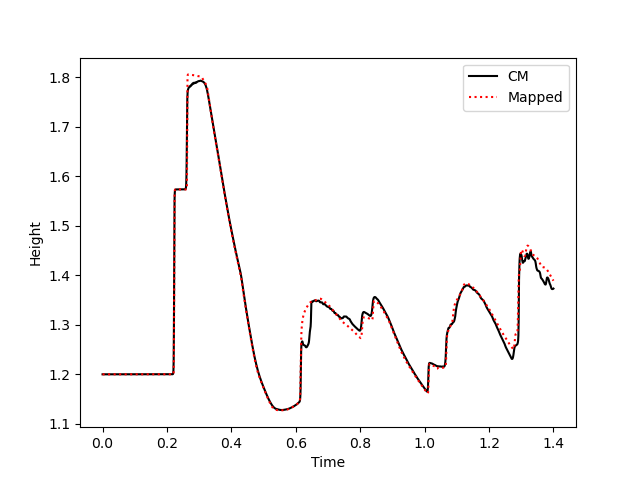}
    \caption{Time profile of Gauge 1 $(0.25,0.6)$.}
    \label{fig:ottpV}
\end{subfigure}
\begin{subfigure}[b]{0.5\textwidth}
\centering
    \includegraphics[scale=0.35]{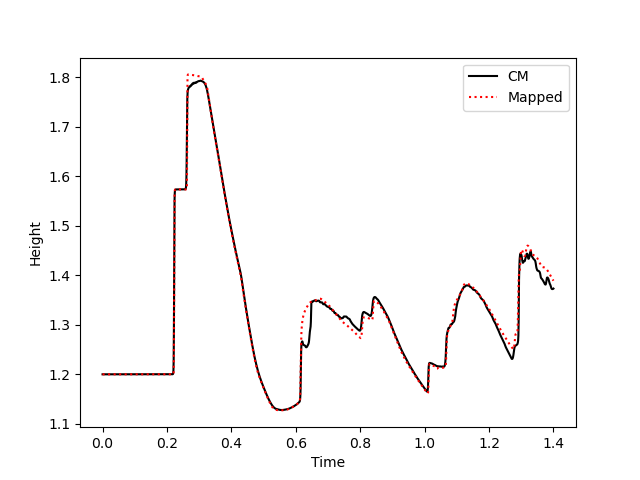}
    \caption{Time profile of Gauge 2 $(0.75,0.6)$.}
    \label{fig:ottpV2}
\end{subfigure}

\begin{subfigure}[b]{0.5\textwidth}
\centering
    \includegraphics[scale=0.35]{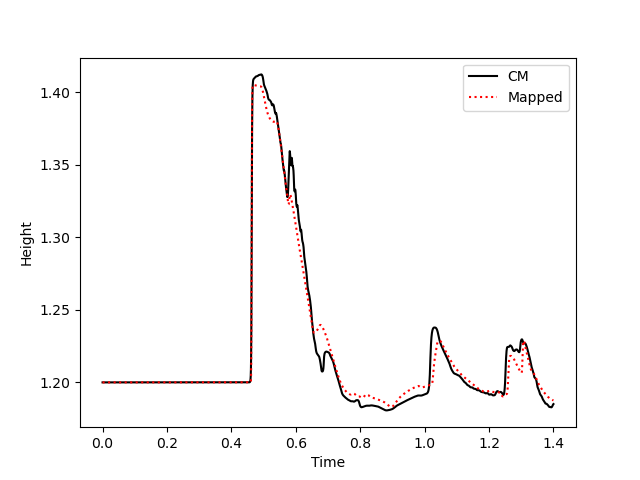}
    \caption{Time profile of Gauge 3 $(0.25,0.3)$.}
    \label{fig:ottpV3}
\end{subfigure}
    \begin{subfigure}[b]{0.5\textwidth}
\centering
    \includegraphics[scale=0.35]{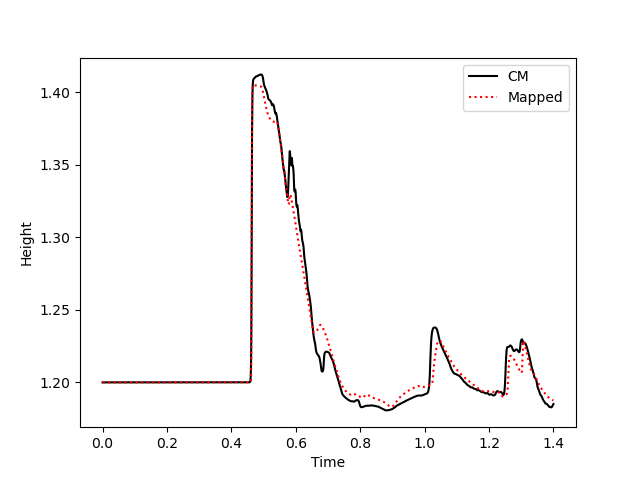}
    \caption{Time profile of Gauge 4 $(0.75,0.3)$.}
    \label{fig:ottpV4}
\end{subfigure}
\caption{Gauge profiles compared with mapped grid results: $900 \times 900$ for CM and $1000 \times 1000$ for mapped grid.}
    \label{fig:otTP}
\end{figure}

We place our gauges at either side of the V-barrier $(0.25,0.3)$, $(0.75,0.3)$, $(0.25,0.6)$, $(0.75,0.6)$ in order to test for symmetry in the results. Indeed we do find symmetry as can be seen in the idential plots of the gauge results in \cref{fig:otTP}.

\subsection{Comparison to mapped grid}
\begin{figure}
    \centering
    \includegraphics[scale=0.5]{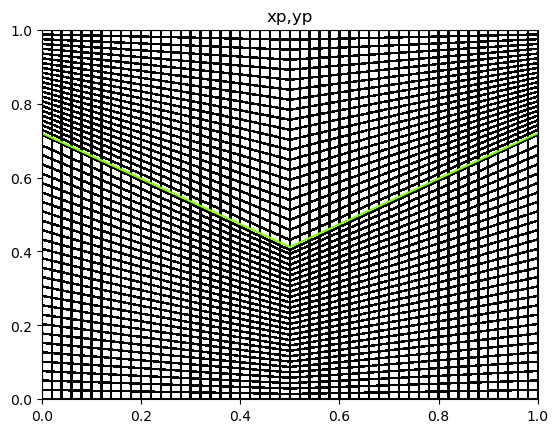}
    \caption{Mapped grid for the V barrier. Coarsened to $50 \times 50$ to highlight mapping. The lime green line represents the edge where the zero width barrier is located. Along this edge we do wave redistribution in the computationally $y$-direction.}
    \label{fig:mapgrid}
\end{figure}

The mapped grid for the V barrier is shown in \cref{fig:mapgrid}. Here we transform the computational uniform grid into a chevron grid akin to what is done in \cite{berger2012simplified}, with the following mapping $f$:
\begin{align}
    f_V(x,y) & = (x, \mu_L (y)), \\
    \mu_V(y) & =\begin{cases} \frac{L_1(x)}{y^*} y  &\mbox{if } (x,y) \in [0,0.5]\times [0,y^*] \\
\frac{1-L_1(x)}{1-y^*}(y-1) +1 & \mbox{if } (x,y) \in [0,0.5]\times [y^*,1] \\
\frac{L_2(x)}{y^*} y  &\mbox{if } (x,y) \in [0.5,1]\times [0,y^*] \\
\frac{1-L_2(x)}{1-y^*}(y-1) +1 & \mbox{if } (x,y) \in [0.5,1]\times [y^*,1]
 \end{cases},
\end{align}
where $y^*$ is the computational barrier edge (taken to be $y^* = 0.5(y_1+y_2) = 0.5(y_2+y_3)$: see \cref{fig:2d_setup2}), and $L_1(x)$ is the barrier line from coordinate 1 to 2 and $L_2(x)$ is the barrier line from coordinate 2 to 3.

The mesh size is chosen such that the barrier in the computational domain lies exactly on $y=y^*$. Transforming the grid with this mapping, using the mapped grid finite volume method explained in \cref{WR} and  \cite{leveque2002finite,calhoun2008logically}, and applying wave redistribution (also \cref{WR}) at the computational barrier edge will give us the numerical solution to the exact same problem that we solve in the Cartesian coordinates using CM cut cell method.
\subsubsection{Convergence}
\begin{figure}[h!]
    \centering
    \includegraphics[scale=0.5]{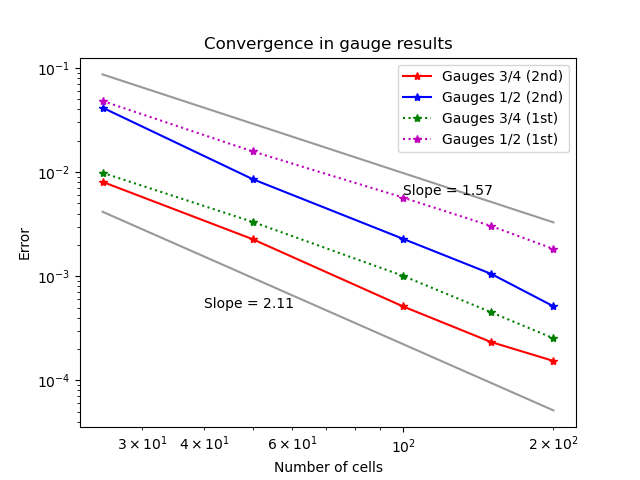}
    \caption{Convergence of gauge profiles.}
    \label{fig:Vconv}
\end{figure}

In \cref{fig:otTP} we plot the gauge results of $1000 \times 1000$ mapped grid V-barrier example and $900 \times 900$ CM example. We do see that the CM results contain more fine movements of the wave and that the mapped grid example produces more smooth wave patterns, implying slight difference in the rate of convergence between the CM method and mapped grid method. Overall, nonetheless, we see convergence as shown in \cref{tab:Vbar_g2}. The order of convergence are around 2 for the second order method and 1.6 for the first order method for both gauge points (\cref{fig:Vconv}).

We note that CM resolves more detail even on a slightly lower resolution than the mapped grid example. This can be seen in the finer details in both the reflected side (where the reflected waves cross in the center) and also the overtopped side, in the appearance of more residual overtopped waves (\cref{fig:ot14_fig}).
\begin{table}[h!]
\centering
\begin{tabular}{|c|c|c|c|c|c|}
\hline
~ & ~ & \multicolumn{2}{c|}{$L_1$ Error (1st)} & \multicolumn{2}{c|}{$L_1$ Error (2nd)}\\
$\Delta x$ & $N_x, N_y$ & Gauge 1/2 & Gauge 3/4 &Gauge 1/2 & Gauge 3/4 \\ \hline
4.e-2   & 25  & 4.83e-2 & 9.86e-3 & 4.16e-2 & 8.07e-3 \\ \hline
2.e-2   & 50  & 1.58e-2 (3.06)& 3.32e-3(2.97) & 8.55e-3(4.87)& 2.27e-3 (3.56) \\ \hline
1.e-2   & 100 & 5.70e-3 (2.77) & 1.01e-3 (3.30)& 2.28e-3 (3.74) & 5.12e-4 (4.42)\\ \hline
0.666e-2 & 150& 3.04e-3 (1.88) & 4.51e-4 (2.23)& 1.05e-3 (2.17)& 2.33e-4 (2.20)\\ \hline
0.5e-2 & 200 & 1.83e-3 (1.66) & 2.53e-4(1.78)& 5.14e-4 (2.05) & 1.54e-4 (1.52)\\ \hline

\end{tabular}
\caption{$L_1$ errors at Gauge 1/2 (0.25,0.6)/(0.75,0.6) and Gauge 3/4 (0.25,0.3)/0.75,0.3).}
\label{tab:Vbar_g2}
\end{table}

\section{Computational Superiority to Refinement using GeoClaw}
\label{sec5}
To highlight the superiority of using the CM method on the zero width barrier, we do a V-barrier simulation on GeoClaw (a geophysical fluid simulator) using same resolution as CM but adaptive (double) refinement at the barrier. In reality the refinement level required at barriers will be greater than 2 as barriers are much skinnier than surrounding bathymetric surfaces. The barrier in these GeoClaw runs is two cells wide, to get down to single cell width in the adaptive mesh refinement.

We show results from using resolution $\Delta x = 1/300$. The lengthened timesteps and number of timesteps show the computational benefit we derive from our proposed method. For $\Delta x =1/300$ we observe that we get 184\% increase in the minimum $\Delta t$ (from 8.6e-06 to 2.4e-05) and 162 \% increase in the average $\Delta t$ and about fivefold decrease in the number of steps taken (from 0.00028 to 0.00074 and 9958 steps to 2037 steps).
\section{Experimental Calculations Comparing the Two Barrier Types}
Here we do final two experiments that will potentially comment on design aspects of a storm barrier in addition to the effectiveness of barriers in general. To compare the effectiveness of the linear barrier versus the V barrier, we add an island on the other side of the barrier and observe flooding at the island center. We will compare the gauge results without a barrier and with a barrier and observe how much protection each barrier type provides.

The boundary conditions are such that there is a wall boundary condition on the side of the wave and extrapolation condition on the side of the island, to allow for wave to exit after hitting the island. Furthermore, for ease of comparison, the V barrier is composed of two copies of the upper half of the linear barrier bent at the midpoint. The initial conditions are $\beta = 1.3$, dam height of $1.5$ and island peak height of 1.3. We provide 3D plots for easier view of the barrier action.

\subsection{Drying and Wetting algorithm}
One side note is the algorithm used for simulating inundation and receding waves on a dry bathymetry. There are many methods that deal with this issue \cite{medeiros2013review}. We employ the relatively simple method described in \cite{bi2014finite} where negatively updated state cells (due to the Riemann problem at water-ground interface) are zeroed out. We note that even though we use the augmented robust Riemann solver known to handle wetting problems, large receding waves (transonic rarefactions) still cause problem in the solver \cite{george2008augmented}. However, adding this fix stabilizes the solver.

\subsection{Linear barrier}
\begin{figure}
    \centering
    \begin{subfigure}[b]{0.35\textwidth}
        \centering

    \includegraphics[scale=0.45]{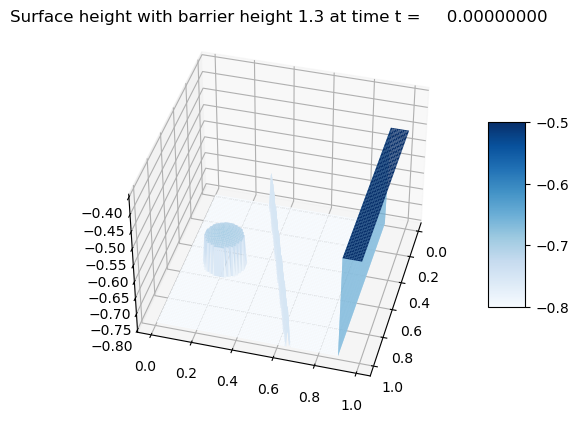}
\caption{Initial}
    \end{subfigure}
    \hspace{2cm}
    \begin{subfigure}[b]{0.35\textwidth}
        \centering

    \includegraphics[scale=0.45]{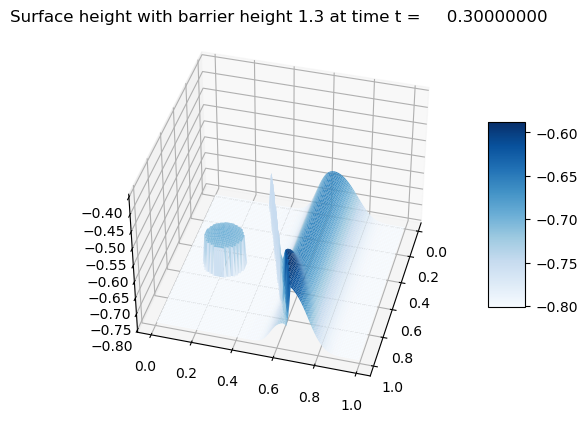}
\caption{First collision against barrier}

    \end{subfigure}
    \vfill
    \begin{subfigure}[b]{0.35\textwidth}
        \centering

    \includegraphics[scale=0.45]{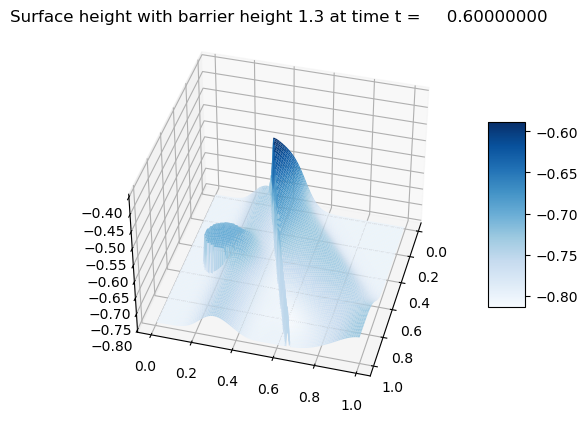}
\caption{First inundation of island}

    \end{subfigure}
    \hspace{2cm}
    \begin{subfigure}[b]{0.35\textwidth}
        \centering

    \includegraphics[scale=0.45]{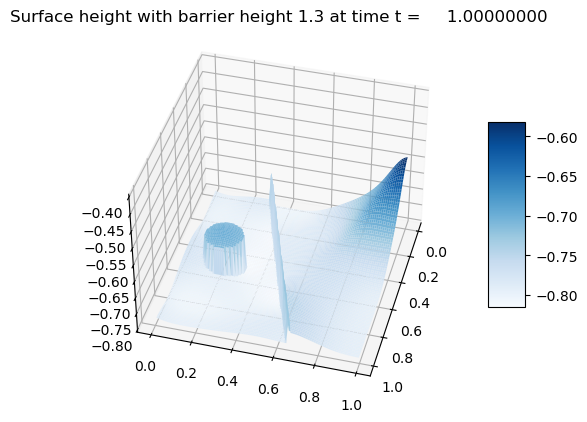}
\caption{Receding and wave leaving domain behind island}

    \end{subfigure}
    \caption{Four time snapshots of barrier action protecting an island. Grid 100 $\times$ 100.}
    \label{fig:L_real}
\end{figure}

\begin{figure}
\begin{subfigure}[b]{0.5\textwidth}
        \centering
\includegraphics[scale=0.5]{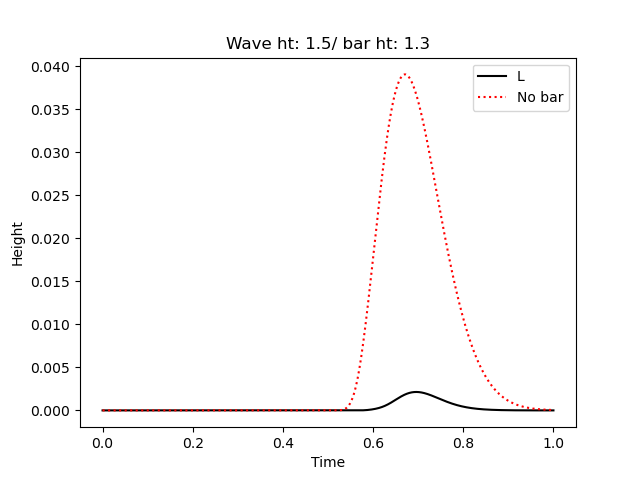}
    \caption{Linear barrier vs no barrier.}
    \label{fig:L_eff}
\end{subfigure}
\begin{subfigure}[b]{0.5\textwidth}
            \centering
\includegraphics[scale=0.5]{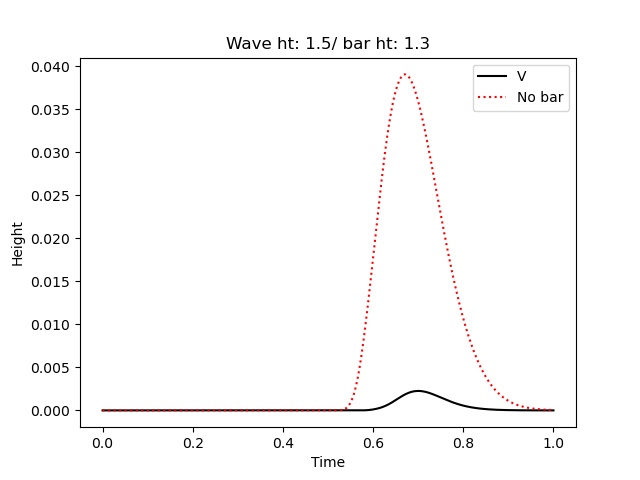}
    \caption{V barrier vs no barrier.}
    \label{fig:V_eff}
\end{subfigure}
\caption{Effectiveness of each barrier measured by gauge results at peak of island.}
\end{figure}

For the linear barrier, we can see from the 3D snapshots \cref{fig:L_real} that the waves glide along the barrier and bound around in the dam break side of the barrier. Inundation starts around $t=0.6.$ The wave is redirected to a direction parallel to the barrier.
We observe that for the linear barrier, inundation is about 0.002. Without the barrier, the inundation becomes around 0.04, implying that the barrier protects the island from about 99.95\% of the peak inundation \cref{fig:L_eff}.

\subsection{V barrier}
\begin{figure}
    \centering
    \begin{subfigure}[b]{0.35\textwidth}
        \centering

    \includegraphics[scale=0.45]{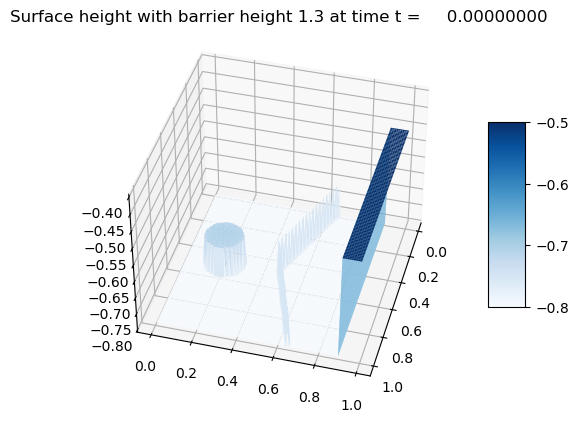}
\caption{Initial}
    \end{subfigure}
    \hspace{2cm}
    \begin{subfigure}[b]{0.35\textwidth}
        \centering

    \includegraphics[scale=0.45]{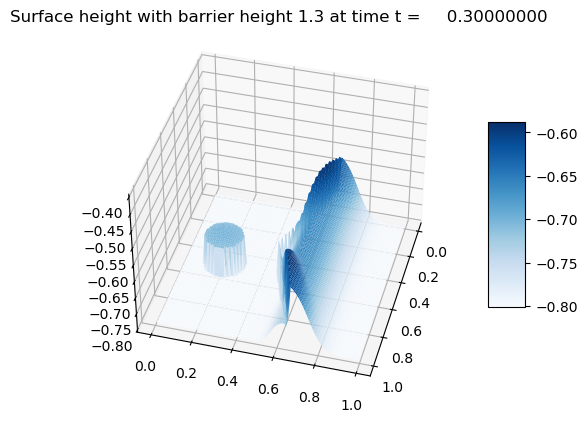}
\caption{First collision against barrier}

    \end{subfigure}
    \vfill
    \begin{subfigure}[b]{0.35\textwidth}
        \centering

    \includegraphics[scale=0.45]{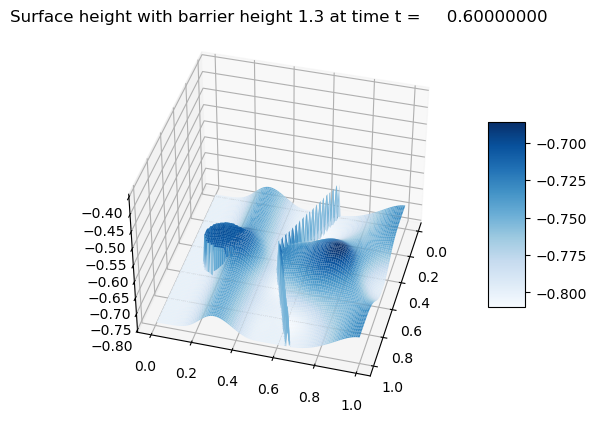}
\caption{First inundation of island}

    \end{subfigure}
    \hspace{2cm}
    \begin{subfigure}[b]{0.35\textwidth}
        \centering

    \includegraphics[scale=0.45]{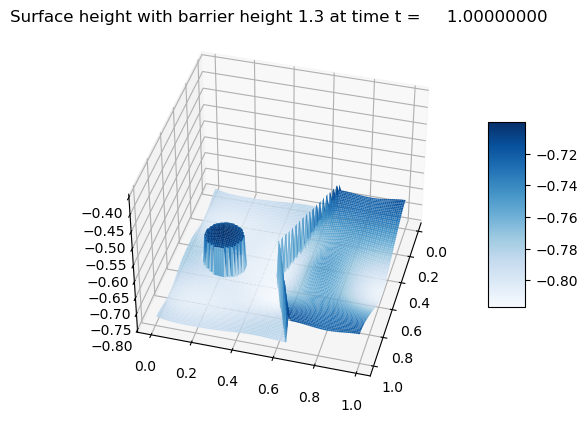}
\caption{Receding and wave leaving domain behind island}

    \end{subfigure}
    \caption{Four time snapshots of barrier action protecting an island. Grid 100 $\times$ 100.}
    \label{fig:V_real}
\end{figure}

\begin{figure}
    \centering
\includegraphics[scale=0.5]{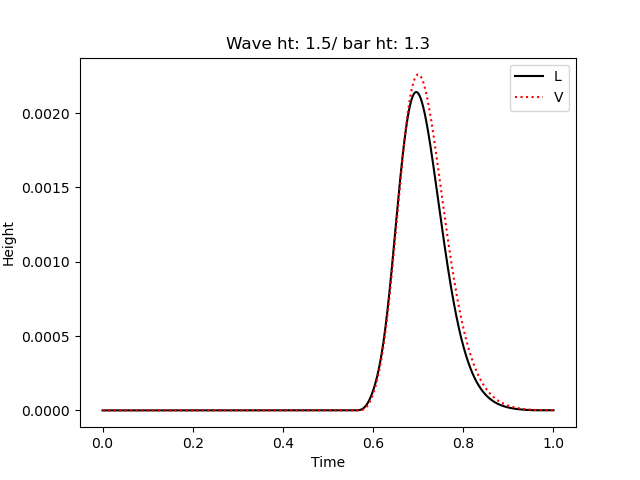}
    \caption{Linear barrier vs V barrier. About the same inundation.}
    \label{fig:LV_Comp}
\end{figure}

We do the same experiment using the V barrier and observe very similar results (\cref{fig:V_real}) when it comes to inundation at the island. The 3D plots show perhaps more controlled behavior when acting against the barrier due to the symmetry (c.f. waves going along barrier and bouncing back etc. for the linear barrier). However, the amount of protection is about the same (~99.95\% of peak flooding \cref{fig:V_eff}), suggesting that it is the height of the barrier that is more determinant in providing protection against flooding (\cref{fig:LV_Comp}). The controlled behavior of the waves in the V barrier could be a merit to consider in barrier design, however, since it directs the incoming wave toward each other toward the center, instead of directing it all in one direction as the linear barrier does. If, for example, there were another island at one end of the linear barrier, the waves could reach there. The V barrier avoids this by aggregating the wave toward the center and evenly distributing out the reflected waves.
\section{Conclusion}
In this paper, we have come up with both first and second order cell merging method handling zero-width barrier simulations. We illustrated the potential for storm simulations by doing a realistic scenario with an island and two types of barrier, the linear and V shaped barrier. We also demonstrated the computational cost saved by doing the proposed method in simulating barriers when compared to resolving the barrier using geophysical fluid simulators such as GeoClaw. We save about 1/5 of the timesteps. Also, we demonstrated that although the design of the barrier does not affect the reduction of inundation at the island (both V and linear protect comparably), the V barrier performs better at directing the waves such that they do not all migrate to one location, where further damage can potentially take place. Finally, having a barrier in our island experiment blocked about 99.95\% of the waves.

\bibliographystyle{siamplain}
\bibliography{bibliography}

\end{document}